%%%

\font\cmrfootnote=cmr10 scaled 750
\font\cmrhalf=cmr10 scaled \magstephalf     
\font\cmrone=cmr10 scaled \magstep1

\font\cmrscfootnote=cmr7 scaled 750
\font\cmrschalf=cmr7 scaled \magstephalf     
\font\cmrscone=cmr7 scaled \magstep1

\font\cmrscscfootnote=cmr5 scaled 750
\font\cmrscschalf=cmr5 scaled \magstephalf     
\font\cmrscscone=cmr5 scaled \magstep1

\font\mitfootnote=cmmi10 scaled 750
\font\mithalf=cmmi10 scaled \magstephalf    
\font\mitone=cmmi10 scaled \magstep1

\font\mitscfootnote=cmmi7 scaled 750
\font\mitschalf=cmmi7 scaled \magstephalf    
\font\mitscone=cmmi7 scaled \magstep1

\font\mitscscfootnote=cmmi5 scaled 750
\font\mitscschalf=cmmi5 scaled \magstephalf    
\font\mitscscone=cmmi5 scaled \magstep1

\font\cmsyfootnote=cmsy10 scaled 750
\font\cmsyhalf=cmsy10 scaled \magstephalf   
\font\cmsyone=cmsy10 scaled \magstep1

\font\cmsyscfootnote=cmsy7 scaled 750
\font\cmsyschalf=cmsy7 scaled \magstephalf   
\font\cmsyscone=cmsy7 scaled \magstep1

\font\cmsyscscfootnote=cmsy5 scaled 750
\font\cmsyscschalf=cmsy5 scaled \magstephalf   
\font\cmsyscscone=cmsy5 scaled \magstep1

\font\cmexfootnote=cmex10 scaled 750
\font\cmexhalf=cmex10 scaled \magstephalf   
\font\cmexone=cmex10 scaled \magstep1

\font\cmexscfootnote=cmex10 scaled 750
\font\cmexschalf=cmex10 scaled \magstephalf   
\font\cmexscone=cmex10 scaled \magstep1

\font\cmexscscfootnote=cmex7 scaled 750
\font\cmexscschalf=cmex7                       
\font\cmexscscone=cmex7 scaled \magstephalf

\def\mathfootnote{\textfont0=\cmrfootnote \textfont1=\mitfootnote \textfont2=\cmsyfootnote \textfont3=\cmexfootnote
        \scriptfont0=\cmrscfootnote \scriptscriptfont0=\cmrscscfootnote \scriptfont1=\mitscfootnote \scriptscriptfont1=\mitscscfootnote
        \scriptfont2=\cmsyscfootnote \scriptscriptfont2=\cmsyscscfootnote \scriptfont3=\cmexscfootnote \scriptscriptfont3=\cmexscscfootnote}
\def\mathhalf{\textfont0=\cmrhalf \textfont1=\mithalf \textfont2=\cmsyhalf \textfont3=\cmexhalf
        \scriptfont0=\cmrschalf \scriptscriptfont0=\cmrscschalf \scriptfont1=\mitschalf \scriptscriptfont1=\mitscschalf
        \scriptfont2=\cmsyschalf \scriptscriptfont2=\cmsyscschalf \scriptfont3=\cmexschalf \scriptscriptfont3=\cmexscschalf}
\def\mathone{\textfont0=\cmrone \textfont1=\mitone \textfont2=\cmsyone \textfont3=\cmexone
        \scriptfont0=\cmrscone \scriptscriptfont0=\cmrscscone \scriptfont1=\mitscone \scriptscriptfont1=\mitscscone
        \scriptfont2=\cmsyscone \scriptscriptfont2=\cmsyscscone \scriptfont3=\cmexscone \scriptscriptfont3=\cmexscscone}

\let\mf=\mathfootnote

\def\meight{\mathone\scriptstyle} 

%%%

\font\Bbb=msbm10
\font\Bbbfootnote=msbm7 scaled \magstephalf

\def\outlin#1{\hbox{\Bbb #1}}
\def\outlinfootnote#1{\hbox{\Bbbfootnote #1}}

\font\call=cmsy10 

\def\cal{\call}
\font\calfoot=cmsy7

\font\small=cmr5                       
\font\notsosmall=cmr7
\font\cmreight=cmr8
 
\font\bfeight=cmbx8

\font\cmrhalf=cmr10 scaled \magstephalf
       
\font\two=cmr10 scaled \magstep2

\font\sansfoot=cmss8
\font\sans=cmss10

\font\twosans=cmss10 scaled \magstep2

\font\caps=cmcsc10  \def\sc{\caps}

%%%

\def\centereps#1#2#3{\vglue#2\relax\centerline{\hbox to#1%
            {\special{eps:#3.eps x=#1 y=#2}\hfil}}}

\let\mf=\mfootnote

\def\q{\quad}   \def\qq{\qquad}

\def\cut{\hfill\break}  \def\newline{\cut}
\def\h#1{\hbox{#1}}

\def\dis{\displaystyle}

\def\eqa{\eqalign} \def\eqano{\eqalignno}

 \def\eps{\epsilon}

\def\a{\alpha}  
 \def\L{\Lambda} 
\def\vp{\varphi}
  
\def\o{\omega}   \def\O{\Omega}

\def\calB{{\h{\cal B}}} \def\calH{{\h{\cal H}}}
\let\H=\calH

\def\calO{{\h{\cal O}}} \def\calD{{\h{\cal D}}}
 
\def\calL{{\h{\cal L}}} \def\L{{\calL}} \def\Lie{{\calL}}
\def\calI{{\h{\cal I}}} \def\calA{{\h{\cal A}}}

\def\calfootH{{\h{\calfoot H}}}  \def\Hfoot{\calfootH}
  
\def\calfootB{{\h{\calfoot B}}}

\def\CC{{\outlin C}} \def\RR{{\outlin R}} \def\ZZ{{\outlin Z}} \def\NN{{\outlin N}}

\def\CCfoot{{\outlinfootnote C}} \def\RRfoot{{\outlinfootnote R}}

\chardef\dotlessi="10 
\chardef\inodot="10

\def\polishl{\char'40l}   
\def\polishL{\leavemode\setbox0=\hbox{L}\hbox to\wd0{\hss\char'40L}}

\def\vol{\h{\rm Vol}}
\def\diam{\hbox{\rm diam}}

\def\Dom{\h{\rm Dom}}

\def\Aut{\h{\rm Aut}}  \def\aut{\h{\rm aut}}  
  
\def\id{\h{\rm id}}

 \def\Dom{\hbox{\rm Dom}}  \def\Ric{\h{\rm Ric}}

\def\tr{\h{\rm tr}}

\def\osc{\h{\rm osc$\,$}} \def\oscfoot{\h{\notsosmall osc$\,$}}

\def\V{\frac1V}

\def\JJJ{\h{\rm J}}

\def\gHerm{g_{\h{\small Herm}}}

\def\gij{g_{i\bar j}}

\def\AutMJ{\hbox{$\hbox{\rm Aut}(M,\hbox{\rm J})$}}

\def\autMJ{\hbox{$\hbox{\rm aut}(M,\hbox{\rm J})$}}

\def\oFSc{\o_{\hbox{\small FS},c}}

\def\GL2nR{\h{$GL(2n,\RR)$}}
\def\Sp2nR{\h{$Sp(2n,\RR)$}}

  \def\dz{d{\bf z}} 
\def\dz#1{dz^{#1}}
\def\dzb#1{d\overline{z^{#1}}}  
\def\dbz{d\bar z}

\def\dzidzjb{dz^i\w\dbz^j}

\def\ddt{\frac d{dt}}

\def\part#1{\frac{\partial#1}{\partial t}}

\def\frac#1#2{{{#1}\over{#2}}}

\def\all{\forall}

\def\sm{\setminus}

\def\precpt
{\hbox{$\mskip3mu\mathhalf\subset\raise0.92pt\hbox{$\mskip-10mu\!\!\!\!
\mathfootnote\subset$}\mskip5mu$}}

\def\supsetnoteq{\hbox{$\mskip3mu\supset\raise-5.97pt
\hbox{$\mskip-10mu\!\!\!\scriptstyle\not=$}\mskip8mu$}}

\def\subsetnoteq{\hbox{$\mskip3mu\subset\raise-5.97pt
\hbox{$\mskip-11mu\!\!\scriptstyle\not=$}\mskip8mu$}}

\def\sseq{\subseteq}

\def\isom{\cong}

\def\*{\star}

\def\w{\wedge}

\def\D{\Delta}
 
\def\dbar{\bar\partial}
\def\del{\partial} \def\pa{\partial}
\def\ddbar{\partial\dbar}

\def\intM{\int_M}  \def\intm{\int_M}
\def\V{\frac1V}
\def\VintM{\V\int_M}

 \def\ha{\frac12}

\def\i{\sqrt{-1}}

\def\ra{\rightarrow}

\def\Ra{\Rightarrow}

\def\arrow#1{\hbox to #1pt{\rightarrowfill}}

\def\thhnotsosmall#1{${\hbox{#1}}^{\hbox{\small th}}$}

\def\nd{2$^{\hbox{\small nd}}$ }

\def\MA{Monge-Amp\`ere }

\def\K{K\"ahler }

\def\KE{K\"ahler-Einstein }  \def\KEno{K\"ahler-Einstein}
\def\KR{K\"ahler-Ricci }

\def\Kolodziej{Ko\polishl{}odziej}

\def\po1{partition of unity }

\def\proof{\hglue-\parindent{\it Proof. }}

\def\Cinf{\h{\cal C}^\infty}

\def\Linf{L^\infty}

\def\Loneloc{L_{\h{\small loc}}^1}

\def\CinfM{\Cinf(M)}

%%%

\def\strutdepth{\dp\strutbox}
\def\specialstar{\vtop to \strutdepth{
    \baselineskip\strutdepth
    \vss\llap{$\star$\ \ \ \ \ \ \ \ \  }\null}}
\def\marginalstar{\strut\vadjust{\kern-\strutdepth\specialstar}}
\def\marginal#1{\strut\vadjust{\kern-\strutdepth
    {\vtop to \strutdepth{
    \baselineskip\strutdepth
    \vss\llap{{ \small #1 }}\null} 
    }}
    }

\newcount\subsectionitemnumber
\def\clearsubsectionitemnumber{\subsectionitemnumber=0\relax}

\newcount\subsubsubsectionnumber
\def\clearsubsubsubsectionnumber{\subsubsubsectionnumber=0\relax}
\def\subsubsubsection#1{
\bigskip\noindent
\global\advance\subsubsubsectionnumber by 1%
{\rm
 \the\subsectionnumber.\the\subsubsectionnumber.\the\subsubsubsectionnumber}
{
#1.}
}

\newcount\subsubsectionnumber
\def\clearsubsubsectionnumber{\subsubsectionnumber=0\relax}
\def\subsubsection#1{
\clearsubsubsubsectionnumber
\bigskip\noindent
\global\advance\subsubsectionnumber by 1%
{%
\it \the\subsectionnumber.\the\subsubsectionnumber}
{
\it #1.}
}

\newcount\subsectionnumber
\def\clearsubsectionnumber{\subsectionnumber=0\relax}
\def\subsection#1{
\clearsubsectionitemnumber
\clearsubsubsectionnumber
\medskip\medskip\smallskip\noindent \global\advance\subsectionnumber by 1%
{%
\bf \the\subsectionnumber} 
{
\bf #1.}
}
\newcount\sectionnumber
\def\clearsectionnumber{\sectionnumber=0\relax}
\def\section#1{
\clearsubsectionnumber
\bigskip\bigskip\noindent \global\advance\sectionnumber by 1%
{%
\two
\the\sectionnumber} 
{
\two #1.}
}

\def\subsubsubsectionno#1{
\bigskip\noindent
\global\advance\subsubsubsectionnumber by 1%
{\rm
 \the\subsectionnumber.\the\subsubsectionnumber.\the\subsubsubsectionnumber}
{
#1}
}
\def\subsubsectionno#1{
\clearsubsubsubsectionnumber
\bigskip\noindent
\global\advance\subsubsectionnumber by 1%
{%
\it \the\subsectionnumber.\the\subsubsectionnumber}
{
\it #1}
}

\def\subsectionno#1{
\clearsubsectionitemnumber
\clearsubsubsectionnumber
\medskip\medskip\smallskip\noindent \global\advance\subsectionnumber by 1%
{%
\bf \the\subsectionnumber} 
{
\bf #1}
}
\def\sectionno#1{
\clearsubsectionnumber
\bigskip\bigskip\noindent \global\advance\sectionnumber by 1%
{%
\two
\the\sectionnumber} 
{
\two #1}
}

\clearsectionnumber
\clearsubsectionnumber
\clearsubsubsectionnumber
\clearsubsubsubsectionnumber
\clearsubsectionitemnumber

\newcount\contentsitemnumber
\def\clearcontentsitemnumber{\contentsitemnumber=0\relax}

\newcount\contentssubitemnumber
\def\clearcontentssubitemnumber{\contentssubitemnumber=0\relax}

\def\contentsitem#1#2{
\clearcontentssubitemnumber
\noindent
\global\advance\contentsitemnumber by 1%
\item{\rm%
\the\contentsitemnumber.}
{%
#1\dotfill #2}
}
\def\contentssubitem#1#2{
\global\advance\contentssubitemnumber by 1%
\itemitem{\rm%
\the\contentsitemnumber.\the\contentssubitemnumber}
{%
#1\dotfill #2}
}

\clearcontentsitemnumber
\clearcontentssubitemnumber

\def\contentsitemno#1#2{
\noindent{#1\dotfill #2} }

\newcount\itemnumber
\def\clearitemnumber{\itemnumber=0\relax}
\def\c#1{ {\noindent\bf \the\itemnumber.} p. $#1$ \global\advance\itemnumber by 1}
\def\cn{ {\noindent\bf \the\itemnumber.} \global\advance\itemnumber by 1}

\clearitemnumber

\def\subsectionno#1{
\medskip\medskip\smallskip\noindent%
{%
\bf #1.}
}

\def\last#1{\advance\eqncount by -#1(\the\eqncount)\advance\eqncount by #1}
\def\llast{\advance\eqncount by -1(\the\eqncount)\advance\eqncount by 1}
\def\lllast{\advance\eqncount by -2(\the\eqncount)\advance\eqncount by 2}
\def\llllast{\advance\eqncount by -3(\the\eqncount)\advance\eqncount by 3}

\newcount\notenumber
\def\clearnotenumber{\notenumber=0\relax}
\def\note#1{\global\advance\notenumber by 1\footnote{${}^{\the\notenumber}$}%
{\lineskip0pt\notsosmall #1}}
\def\notewithcomma#1{\advance\notenumber by
1\footnote{${}^{\the\notenumber}$}
  {\lineskip0pt\notsosmall #1}}

\clearnotenumber

\def\putnumber{%
\global\advance\subsectionitemnumber by 1{\the\subsectionnumber}.{\the\subsectionitemnumber}}
\def\numbering{{{\the\subsectionnumber}.{\the\subsectionitemnumber}}}

\def\Abstract#1{
{\narrower\bigskip\bigskip
\noindent {{\bf Abstract.\ \ } #1}

}
}

\def\FThm#1{\bigskip
{\noindent%
{\bf {T}heorem }
{\bf \putnumber.} {\it #1}\bigskip 
}}

\def\FCor#1{\bigskip
{\noindent%
{\bf {C}orollary }
{\bf \putnumber.} {\it #1}\bigskip
}}
\def\FLem#1{\bigskip
{\noindent%
{\bf {L}emma }
{\bf \putnumber.} {\it #1}\bigskip
}}

\def\FProp#1{\bigskip
{\noindent%
{\bf {P}roposition }
{\bf \putnumber.} {\it #1}\bigskip
}}

\def\FRemm#1{\medskip
{\noindent%
{\bf {R}emark }
{\bf \putnumber.}  {#1}\medskip
}}

\def\FConj#1{\bigskip
{\noindent%
{\bf {C}onjecture }
{\bf \putnumber.} {\it #1}\bigskip
}}

\def\FDef#1{\bigskip
{\noindent%
{\bf {D}efinition }
{\bf \putnumber.} {\it #1}\bigskip
}}

\def\ref#1{{\bf[}{\sans #1}{\bf]}}

\def\reffoot#1{{\bfeight[}{\sansfoot #1}{\bfeight]}}

\def\opcit{\underbar{\phantom{aaaaa}}}

\def\sm{\smallskip}

\def\boxit#1{\vbox{\hrule\hbox{\vrule\kern3pt\vbox{\kern3pt#1\vglue3pt
\kern3pt}\kern3pt\vrule}\hrule}}

\long\def\frame#1#2#3#4{\hbox{\vbox{\hrule height#1pt
 \hbox{\vrule width#1pt\kern #2pt
 \vbox{\kern #2pt
 \vbox{\hsize #3\noindent #4}
\kern#2pt}
 \kern#2pt\vrule width #1pt}
 \hrule height0pt depth#1pt}}}

\def\kou{\frame{.3}{.5}{1pt}{\phantom{a}}}

\def\help{\ifmmode\aftergroup\noindent\quad\else\quad\fi}
\def\helpp{\ifmmode\aftergroup\noindent\hfill\else\hfill\fi}
\def\done{\helpp\kou\medskip}

\def\smallblackbox{\vrule height.6ex width .5ex depth -.1ex}
\def\boxseparation{\hfil\smallblackbox$\q$\smallblackbox$\q$\smallblackbox\hfil}

%%%

%%

\def\Ho{{\calH_{\o}}}

\def\on{\o^n}
\def\onminus1{\o^{n-1}/(n-1)!}
\def\ovp{\o_{\vp}}

\def\vpt{\vp_{\!t}}

\def\ovpt{\o_{\vpt}}

\def\ovpn{\o_{\vp}^n}
\def\ovptn{\o_{\vpt}^n}

%%%

\font\call=cmsy10   \font\Bbb=msbm10
\def\outlin#1{\hbox{\Bbb #1}}
\def\H{\hbox{\call H}}  \def\Ho{\hbox{\call H}_\omega}
 
\def\Cinf{{\call C}^\infty}  \def\RR{{\outlin R}}
\def\frac#1#2{{{#1}\over{#2}}}
\newcount \eqncount
\def \eqnno{\global \advance \eqncount by 1 \futurelet \nexttok \parsenexttok}
\def \eqn{\global \advance \eqncount by 1 \eqno\futurelet \nexttok \parsenexttok}
\def \eqnd{\global\advance \eqncount by 1 \futurelet\nexttok\parsenexttokd}
\def \parsenexttok{\ifx \nexttok $\Nomark\else\expandafter \Mark\fi}
\def \parsenexttokd{\ifx \nexttok \hfil\Nomark\else\expandafter \Mark\fi}
\def \Nomark {(\the \eqncount)}\def \Mark #1{\xdef #1{(\the \eqncount)}#1}

%%%

\let\exa\expandafter
\catcode`\@=11 
\def\CrossWord#1#2#3{%
\def\@x@{}\def\@y@{#2}
\ifx\@x@\@y@ \def\@z@{#3}\else \def\@z@{#2}\fi
\exa\edef\csname cw@#1\endcsname{\@z@}}
\openin15=\jobname.ref
\ifeof15 \immediate\write16{No file \jobname.ref}%
   \else \input \jobname.ref \fi 
\closein15
\newwrite\refout
\openout\refout=\jobname.ref 
\def\warning#1{\immediate\write16{1.\the\inputlineno -- warning --#1}}
\def\Ref#1{%
\exa \ifx \csname cw@#1\endcsname \relax
\warning{\string\Ref\string{\string#1\string}?}%
    \hbox{$???$}%
\else \csname cw@#1\endcsname \fi}
\def\Tag#1#2{\begingroup
\edef\head{\string\CrossWord{#1}{#2}}%
\def\writeref{\write\refout}%
\exa \exa \exa
\writeref \exa{\head{\the\pageno}}%
\endgroup}

\catcode`\@=12

\def\Tagg#1#2{\Tag{#1}{#2}}  

\def\TaggThm#1{\Tagg{#1}{Theorem~\numbering}}
\def\TaggLemm#1{\Tagg{#1}{Lemma~\numbering}}
\def\TaggRemm#1{\Tagg{#1}{Remark~\numbering}}
\def\TaggCor#1{\Tagg{#1}{Corollary~\numbering}}
\def\TaggProp#1{\Tagg{#1}{Proposition~\numbering}}
\def\TaggSection#1{\Tagg{#1}{Section~\the\subsectionnumber}}
\def\TaggS#1{\Tagg{#1}{\S\the\subsectionnumber}}
\def\TaggSubS#1{\Tagg{#1}{\hbox{\S\S\unskip\the\subsectionnumber.\the\subsubsectionnumber}}}
\def\TaggSubsection#1{\Tagg{#1}{Subsection~\the\subsectionnumber.\the\subsubsectionnumber}}
\def\TaggEq#1{\Tagg{#1}{(\the\eqncount)}}
\def\Taggf#1{ \Tagg{#1}{{\bf[}{\sans #1}{\bf]}} }
\def\TaggDef#1{\Tagg{#1}{Definition~\numbering}}

\newcount\pagenotemp
\def\pagenotempshow{\pagenotemp=\pageno\relax \advance\pagenotemp by 1\relax}
\def\clearpagenotemp{\advance\pagenotemp by -1\relax}
\def\TaggPage#1{\pagenotempshow\Tagg{#1}{\the\pageno%
}\clearpagenotemp}
  
\def\TaggConj#1{\Tagg{#1}{Conjecture~\numbering}}

%%%

\def\ricm#1{\hbox{\rm Ric}^{-#1}}
\def\Ric{\hbox{\rm Ric}\,}
\def\Ricnotsosmall{\hbox{\notsosmall Ric}\,}

\def\Hc{\H_{c_1}}
\def\Dc{\calD_{c_1}}

\def\Hcfoot{\Hfoot_{c_1}}

\def\hc#1{\H_{c_1}^{(#1)}}

\def\HcG{\H_{c_1}(G)}

\def\HcGplus{\H^{+}_{c_1}(G)}
\def\Hcplus{\H^{+}_{c_1}}

\def\Hmuc{\H_{\mu c_1}}

\def\O{\Omega}

\def\Rico{\Ric\o}
\def\Ricovp{\Ric\ovp}

\def\Ho{\H_{\Omega}}
\def\Do{\calD_{\Omega}}
\def\ho#1{\H_{\Omega}^{(#1)}}

\def\HoG{\H_{\Omega}(G)}

\def\MA{Monge-Amp\`ere }

\def\Vm{V^{-1}}  \def\V{\frac1V}

\def\fovp{f_{\ovp}}
\def\fo{f_{\o}}
\def\fovpt{f_{\ovpt}}

%%%%%%%

\def\JJJ{\h{\rm J}}
\def\Ric{\hbox{\rm Ric}\,}
\def\Ricnotsosmall{\hbox{\notsosmall Ric}\,}

\def\Hc{\H_{c_1}}
\def\Dc{\calD_{c_1}}
\def\Hcfoot{\Hfoot_{c_1}}
\def\HcG{\H_{c_1}(G)}

\def\HcGplus{\H^{+}_{c_1}(G)}
\def\Hcplus{\H^{+}_{c_1}}

\def\Hmuc{\H_{\mu c_1}}

\def\O{\Omega}

\def\Rico{\Ric\o}
\def\Ricovp{\Ric\ovp}

\def\HO{\H_{\Omega}}
\def\Ho{{\calH_{\o}}}
\def\Do{\calD_{\Omega}}
\def\HoG{\H_{\Omega}(G)}

\def\MA{Monge-Amp\`ere }

\def\Vm{V^{-1}}  \def\V{\frac1V}

\def\fovp{f_{\ovp}}
\def\fo{f_{\o}}

%%%

\def\KR{K\"ahler-Ricci }

\magnification=1100
\hoffset1.1cm
\voffset1.3cm
\hsize5.09in
\vsize6.91in

\headline={\ifnum\pageno>1{\ifodd\pageno \oddheadline\else\evenheadline\fi}\fi}
\def\oddheadline{\centerline{\caps Some discretizations of geometric evolution equations}}
\def\evenheadline{\centerline{\caps Y. A. Rubinstein}}

\overfullrule0pt
\parindent12pt
\vglue0.6cm

\noindent
\centerline{\twosans Some discretizations of geometric evolution equations}
\medskip
\centerline{\twosans and the Ricci iteration on the space of K\"ahler metrics, I}
\smallskip

\font\tteight=cmtt8 

\bigskip
\centerline{\sans Yanir A. Rubinstein%
\footnote{$^*$}{\cmreight Massachusetts Institute of Technology. 
Email: {\tteight yanir@member.ams.org}}}

\def\centereps#1#2#3{\vglue#2\relax\centerline{\hbox to#1%
            {\special{eps:#3.eps x=#1 y=#2}\hfil}}}

\font\cmreight=cmr8

\vglue0.03in

\footnote{}{\hglue-\parindent\cmreight September \thhnotsosmall{6}, 2007. 
\hfill\break 
\vglue-0.5cm
\hglue-\parindent\cmreight
\vglue-0.5cm 
\hglue-\parindent\cmreight
32Q20, % 
39A12, %
53C25, %   
58E11.} %  
%       % 
    
\vglue-0.12in

\bigskip
\Abstract{
In this article and in its sequel
we propose the study of certain discretizations of geometric evolution equations
as an approach to the study of the existence problem of some elliptic partial differential equations
of a geometric nature as well as a means to obtain interesting dynamics on certain infinite-dimensional
spaces.
We illustrate the fruitfulness of this approach in the context of the Ricci flow, as well as another
flow, 
in K\"ahler geometry. We introduce and study dynamical systems related to 
the Ricci operator on the space of \K metrics that arise as discretizations of these flows.
We pose some problems regarding their dynamics.
We point out a number of applications to well-studied objects in \K and conformal geometry such as
constant scalar curvature metrics, \KR solitons, Nadel-type multiplier ideal sheaves, balanced metrics, the Moser-Trudinger-Onofri inequality, energy functionals and the geometry and structure of the space of \K metrics.
E.g., we obtain a new sharp inequality strengthening the classical Moser-Trudinger-Onofri inequality on the two-sphere.
}

\bigskip

\medskip
\subsectionnumber=0\relax
{\narrower
\centerline{\sc Contents}
\bigskip

\contentsitem
{Introduction}
{\Ref{IntroductionPage}}

\contentsitem
{Constructing canonical metrics in \K geometry}
{\Ref{CanonicalMetricsPage}} 

\contentsitem
{The Ricci iteration}
{\Ref{RicciIterationPage}} 

\contentsitem
{Some energy functionals on the space of \K metrics}
{\Ref{EnergyFunctionalsPage}} 

\contentsitem
{The Ricci iteration for negative and zero first Chern class}
{\Ref{RicciIterationNegativeZeroPage}} 

\contentsitem
{The Ricci iteration for positive first Chern class}
{\Ref{RicciIterationPositivePage}}

\contentsitem
{The \KR flow and the Ricci iteration for a general \K class}
{\Ref{RicciIterationGeneralPage}}

\contentsitem
{Another flow and the inverse Ricci operator for a general \K class}
{\Ref{AnotherFlowPage}}

\contentsitem
{The twisted Ricci iteration and a twisted inverse Ricci operator}
{\Ref{RicciIterationTwistedPage}}

\contentsitem
{Some applications}
{\Ref{ApplicationsPage}}

\contentssubitem
{The Moser-Trudinger-Onofri inequality on the Riemann sphere\newline and its higher-dimensional analogues}
{\Ref{MTOPage}}

\contentssubitem
{An analytic characterization of \KE manifolds and an\newline analytic criterion for almost-\KE manifolds}
{\Ref{AnalyticCharacterizationPage}}

\contentssubitem
{A new Moser-Trudinger-Onofri inequality on the Riemann sphere\newline and a family of energy functionals}
{\Ref{ImprovedMTOEnergyFunctionalsPage}}

\contentssubitem
{Construction of Nadel-type obstruction sheaves}
{\Ref{NadelSheavesPage}}

\contentssubitem
{Relation to balanced metrics}
{\Ref{BalancedMetricsPage}}

\contentssubitem
{A question of Nadel}
{\Ref{NadelQuestionPage}}

\contentssubitem
{The Ricci index and a canonical nested structure on the space of\newline \K metrics}
{\Ref{RicciIndexPage}}

\contentsitemno
{Bibliography}{\Ref{BibliographyPage}}

}
\bigskip

\subsection{Introduction}{}
\TaggSection{SectionIntroduction}%
{\TaggPage{IntroductionPage}}%
Our main purpose in this article and in its sequel \ref{R5} is to propose the systematic use of certain discretizations 
of geometric evolution equations as an approach to the study of the existence problem of certain elliptic partial 
differential equations of a geometric nature
as well as a means to obtain interesting dynamics on certain infinite-dimensional
spaces. We illustrate the fruitfulness of this approach in the context
of the Ricci flow, as well as another flow, in \K geometry. We describe how this approach 
gives a new method for the construction of canonical \K metrics. We also introduce
a number of canonical dynamical systems on the space of \K metrics that
we believe merit further study. Some of the results and constructions described here were announced previously \ref{R3}.

Given an elliptic partial differential equation, several classical methods are available
to approach the problem of existence of solutions. In essence, standard elliptic theory 
reduces the existence problem to the demonstration of certain a priori estimates for
solutions. The main difficulty lies therefore in devising methods to obtain 
these estimates.

One common method, that goes back at least to Bernstein and Poincar\'e, 
is the continuity method. In this approach one continuously  
deforms the given elliptic operator to
another (oftentimes in a linear fashion), for which the existence problem is known to have solutions. Ellipticity
provides for existence of solutions for small perturbations of this easier problem. In order to
prove existence for the whole deformation path one then seeks to establish
a priori estimates, uniform along the deformation, for solutions of the family of elliptic problems.

Another approach, drawing some of its motivation from Physics, is the heat flow method,
going back to Fourier.
Here the idea is study a deformation of the elliptic problem according to a parabolic
heat equation whose equilibrium state is precisely a solution to the original elliptic equation.
Much of the standard elliptic theory has a parabolic counterpart. First, one makes use of the latter
in order to establish short-time existence. Long-time existence and convergence then
hinge upon establishing a priori estimates, as before.

A third approach, going back, among others, to Euler and Cauchy, is the discretization method, that
can be considered as a blend of the two above.
Here the idea  
is to replace an evolution equation (or ``flow") by a countable set of elliptic equations
that arise by repeatedly solving a difference equation corresponding to discretizing the flow equation
in the time variable.
This approach provides common and elementary numerical algorithms, the Euler
method and its variants, and is widely used in the ``real world", for example to numerically integrate
differential equations.

In this article we wish to explore this third approach in the context of certain geometric evolution equations.
To the best of our knowledge, it seems that it has not been used before 
in a systematic manner in this context.

We would like to emphasize that when a particular elliptic equation has a solution one morally expects
all three methods to converge towards such a solution. Therefore one should not take as a surprise
the fact that the discretization method converges in some of the cases we consider. The crux is thus 
not the convergence itself but rather the new point of view and insights that this method provides; 
both to the study of the original elliptic problem as well as to the understanding of the evolution equation, the
continuity method and the relation between the two. In addition, one may obtain in this way
non-trivial canonical discrete dynamical systems on infinite-dimensional spaces that may be of
some interest in their own right.

Let us consider as a simple illustration the Laplace equation on a bounded 
smooth domain $\O$ in $\RR^n$. The elliptic problem is then to find a function $u$ satisfying
$$
\eqa{
\D u & = 0,\q \h{on\ } \O,\cr
u & = \psi,\q \h{on\ } \pa\O.
}\eqn\LaplaceEq
$$
Consider then the difference equations
$$
\eqa{
u_k-u_{k-1} & = \D u_k,\q \h{on\ } \O,\cr
u_k & = \psi,\q \h{on\ } \pa\O,\cr
u_0 & = u,
}
$$
where $u$ is any smooth function that agrees with the smooth function $\psi$ on the boundary. Thus, one may
write
$$
u_k=(-\D+1)^{-1}\circ\cdots\circ(-\D+1)^{-1} u.
$$
One may then readily show that the sequence $\{u_k\}_{k\ge0}$ exists for each $k\in\NN$ 
and converges exponentially fast in $k$ to the unique solution of \LaplaceEq.

To give further intuition as to why this method works we consider the following finite-dimensional 
problem:
Given a positive semi-definite matrix $A$ and a vector $v$, find the projection of $v$
onto the zero eigenspace of $A$. One possible solution is to consider the sequence of 
vectors defined iteratively by
$$\eqa{
v_k & =(A+I)^{-1}v_{k-1},\cr
v_0 & = v.
}$$
Then $\lim_{k\ra\infty} v_k$ exists and is the required projection.
This algorithm is nothing but the discretization of the flow
$$\eqa{
\frac{d v(t)}{dt} & = -Av(t),\cr
v(0) & = v.
}$$
Discretizations corresponding to different time steps will produce equivalent dynamical systems
$$\eqa{
v_k & =(\tau A+I)^{-1}v_{k-1},\cr
v_0 & = v,
}$$
whose convergence is faster the larger the time-step $\tau\in(0,\infty)$, with $\tau$ and
the first non-zero eigenvalue of $A$ controlling the exponential factor of the speed of convergence.

\subsection{Constructing canonical metrics in \K geometry} 
\TaggSection{CanonicalMetricsSection}%
{\TaggPage{CanonicalMetricsPage}}% 
In this article we wish to apply the method described in the previous section
towards the study of canonical \K metrics and the space of \K metrics. In this section we very briefly describe 
the problem and some background. We refer to \ref{A3,Bes,F1,Si,T6} for more background.

The search for a canonical metric representative of 
a fixed \K class has been at the heart of \K geometry since its birth. Indeed, in his visionary article
\K defined the eponymous manifold motivated by the fact that in this setting Einstein's equation
simplifies considerably and reduces to a second order partial differential equation for a single function
\ref{K}, namely, the local potential $u$ which represents the \K form $\o$ 
on the open domain $U$ via $\o|_U=\i\ddbar u$ must satisfy
$$
\det \Bigg[\frac{\pa^2 u}{\pa z^i \pa{\bar z}^j}\Bigg]= e^{-\mu u},\q \h{on\ } U,
$$
where $\mu$ is the Einstein constant.
Two decades later, following Chern's fundamental work on characteristic classes, 
Calabi introduced the concept of the space of \K metrics in a fixed 
cohomology class and
formulated the problem on a compact closed manifold as an equation for a global
smooth function (\K potential) $\vp$, 
$$
\ovpn=\on e^{\fo-\mu\vp},
$$
where $\fo$ satisfies $\i\ddbar\fo=\Ric\o-\mu\o$. This showed that a necessary condition
for the existence of solutions is that the first Chern class be definite or zero.
Calabi proposed that this
equation should always admit a unique solution in each \K class when $\mu=0$. In addition he suggested 
the study of a more general notion, that of an extremal metric \ref{C1,C2}. 
Since then much progress has been made towards understanding when such metrics exist. Regarding \KE metrics, 
the most general result
in this direction is given by the work of Aubin in the negative Ricci curvature case \ref{A}
and by Yau in the case of nonpositive Ricci curvature which provided a solution to Calabi's conjecture \ref{C1,Y}. 
Following this much work has gone into understanding the positive case, notably by
Tian who provided a complete solution for complex surfaces, in addition to establishing 
an analytic characterization of \KE manifolds and a theory of stability \ref{T1,T2}. 
For general extremal metrics however a general existence theory is not presently available although a 
conjectural picture, the so-called Yau-Tian-Donaldson conjecture, suggests that it should be
intimately related with notions of stability in algebraic geometry \ref{Th}.

The principal tool in the study of \KE metrics has been the continuity method, as suggested 
initially by Calabi \ref{C2}, and later studied by Aubin and Yau.
In the remaining case ($\mu>0$) Bando and Mabuchi showed 
that the continuity method will converge to a \KE metric when one exists \ref{BM}.
Another important tool has been the Ricci flow introduced 
by Hamilton \ref{H} in the more general setting of Riemannian manifolds. 
Cao has shown that the continuity method proofs for the cases $\mu\le0$ may be phrased in terms 
of the convergence of the Ricci flow \ref{Ca}. Later, much work has gone into understanding the 
Ricci flow on Fano manifolds and recently Perelman and Tian and Zhu proved that the analogous 
convergence result holds in this case \ref{TZ4}.

The idea that there might be another way of approaching canonical metrics, in the form of a discrete 
infinite-dimensional iterative dynamical system, was suggested by Nadel \ref{N2}.
More recently, Donaldson has proposed a program for the construction of constant scalar curvature \K
metrics on projective manifolds using finite-dimensional iteration schemes and balanced metrics which are in essence computable
\ref{Do3}. 

The main motivation for our work came from trying to approach Nadel's basic problem: Find 
an (infinite-dimensional) iterative dynamical system on the space of \K metrics that 
converges to a \KE metric.
In his note Nadel suggested one such dynamical system on Fano manifolds which, as we explain below, is not 
suited to the problem  (see \Ref{NadelQuestionSubS}). Nevertheless his idea is related
to the right answer, and it is our purpose in this article and its sequel to describe our approach to Nadel's
problem and some of its consequences.

\subsectionno{Setup and notation}
Let $(M,\JJJ,g)$ be a connected compact closed \K manifold of complex dimension $n$
and let $\O\in H^2(M,\RR)\cap H^{1,1}(M,\CC)$ be a \K class
with $d=\del+\dbar$. Define the Laplacian $\D=-\dbar\circ\dbar^\star-\dbar^\star\circ\dbar$ with
respect to a Riemannian metric $g$
and assume that $\JJJ$ is compatible with $g$ and parallel
with respect to its Levi-Civita connection. 
Let
$\gHerm=1/\pi\cdot\gij(z)\dz i\otimes\dzb j$ be the associated \K metric, 
that is the induced Hermitian metric on $(T^{1,0}M,\JJJ)$, and let
$\o:=\o_g=\i/2\pi\cdot\gij(z)\dzidzjb$ denote its corresponding \K form,
a closed positive $(1,1)$-form on $(M,\h{\rm J})$
such that $\gHerm=\ha g-\frac\i2\o$. Similarly denote by $g_\o$ the Riemannian
metric induced from $\o$ by $g_\o(\cdot,\cdot)=\o(\cdot,\JJJ\,\cdot)$. 

For any \K
form we let $\Ric(\omega)=-\i/2\pi\cdot\ddbar\log\det(\gij)$ denote the
Ricci form of $\o$. It is well-defined globally and represents
the first Chern class $c_1:=c_1(T^{1,0}M,\h{\rm J})\in H^2(M,\ZZ)$. 
Alternatively it may be viewed as minus the curvature form of the canonical line bundle 
$K_M$, the top exterior product of the holomorphic cotangent bundle $T^{1,0\,\star}M$.
One calls $\o$ \KE if $\Ric\o=a\o$ for some real $a$. The trace of the Ricci form
with respect to $\o$ is called the scalar curvature and is denoted by $s(\o)$. The
average of the scalar curvature is denoted by $s_0$ and does not depend on the choice
of $\o\in\HO$. Nor does the volume $V=\O^n([M])$.

Let $H_g$ denote the Hodge projection operator from
the space of closed forms onto the kernel of $\D$. 
Denote by $\Do$ the space of all closed $(1,1)$-forms cohomologous to $\O$, by
$\HO$ the subspace of \K forms, and by $\HO^+$ the 
subspace of \K forms whose Ricci curvature is positive (nonempty if and only if $c_1>0$).

For a \K form $\o$ with $[\o]=\O$ we will consider the space of smooth strictly $\o$-plurisubharmonic 
functions (\K potentials)
$$
\calH_\o=\{\vp\in\CinfM\,:\, \ovp:=\o+\i\ddbar\vp >0\}.
$$
Let $\Aut(M,\JJJ)$ denote the complex Lie group of automorphisms (biholomorphisms)
of $(M,\JJJ)$
and denote by  $\aut(M,\JJJ)$ its Lie algebra of infinitesimal automorphisms composed
of real vector fields $X$ satisfying $\L_X\JJJ=0$. 
Let $G$ be 
any compact real Lie subgroup of $\Aut(M,\JJJ)$, and 
let $\Aut(M,\JJJ)_0$ denote the identity component
of $\Aut(M,\JJJ)$.
We denote by $\HO(G)\sseq\HO$ and $\Ho(G)\sseq\Ho$ the corresponding subspaces of $G$-invariant elements.

\subsection{The Ricci iteration}{} 
\TaggSection{SectionRicciIteration}%
{\TaggPage{RicciIterationPage}}% 
In this section we introduce the Ricci iteration and describe some of its elementary properties.

Hamilton's Ricci flow on a \K manifold of definite or zero first Chern class is defined as the
the set $\{\o(t)\}_{t\in\RRfoot_+}$ satisfying the evolution equations
$$
\eqa{
\frac{\pa \o(t)}{\pa t} & =-\Ric\o(t)+\mu\o(t),\quad t\in\RR_+,\cr
\o(0) & =\o\in\HO,
}\eqn\RicciFlowEq$$
where $\O$ is a \K class satisfying $\mu\O=c_1$ for some $\mu\in\RR$ (see, e.g., \ref{Cho}).
We will sometimes refer to the equations \RicciFlowEq\ themselves as the Ricci flow.

We introduce the following dynamical system that 
is our main object of study in this article. It is a discrete version of this flow.

\FDef
{
\TaggDef{RicciIterationDef}%
Let $\O$ be a \K class satisfying $\mu\O=c_1$ for some $\mu\in\RR$.
Given a \K form $\o\in\HO$ and a number $\tau>0$ define the time $\tau$ Ricci iteration %
 to be 
the sequence of forms $\{\o_{k\tau}\}_{k\ge0}\sseq\HO$,
satisfying the equations
$$
\eqa{
\o_{k\tau} & =\o_{(k-1)\tau}+\tau\mu\o_{k\tau}-\tau\Ric\o_{k\tau},\quad   k\in\NN,\cr
\o_0 & =\o,
}\eqn\IterationTauEq$$
for each $k\in\NN$ for which a solution exists in $\HO$. 
}

We pose the following elementary conjecture concerning the limiting behavior of the Ricci iteration in
the presence of fixed points. 

\FConj{
\TaggConj{MainConjecture}
Let $(M,\JJJ)$ be a compact \K manifold admitting a \KE metric. 
Let $\O$ be a \K class such that $\mu\O=c_1$ with $\mu\in\RR$.
Then for any $\o\in\HO$ and for any $\tau>0$, the time $\tau$ Ricci iteration exists for all $k\in\NN$ and converges
in the sense of Cheeger-Gromov 
to a \KE metric.
}

Regarding this conjecture we prove in this article the following simple result.%
\note{\cmreight We refer the reader to the sequel \reffoot{R5} where we intend to discuss 
some of the remaining cases.\break 
\vglue-0.5cm
\hglue-\parindent\cmreight These cases that we also find of significant geometric interest require more involved and detailled\break 
\vglue-0.5cm
\hglue-\parindent\cmreight 
analysis. Also we wished to keep this first article's length reasonable as well as not to further delay\break 
\vglue-0.5cm
\hglue-\parindent\cmreight this first part that introduces our main constructions and ideas regarding discrete dynamical systems\break 
\vglue-0.5cm
\hglue-\parindent\cmreight on the space of \K metrics.}
\FProp{
\TaggProp{ConvergenceProp}%
Let $(M,\JJJ)$ be a compact \K manifold admitting a unique \KE metric. 
Let $\O$ be a \K class such that $\mu\O=c_1$ with $\mu\in\RR$.
Then for any $\o\in\HO$ and for any $\tau>1/\mu$ (when $\mu>0$) or $\tau>0$ (when $\mu\le0$), the time $\tau$ Ricci iteration exists for all $k\in\NN$ and converges
to a \KE metric.
}

Now, let the Ricci potential $f:\vp\in\Ho\ra \fovp\in\CinfM$ be the vector field on $\Ho$
satisfying
$$
\i\ddbar\fovp=\Ric\ovp-\ovp,\q \V\intm e^{\fovp}\ovpn=1.\eqn\RicciPotentialEq
$$ 
For each $k$ write 
$$
\o_{k\tau}=\o_{\psi_{k\tau}},\q \h{\ with \ } \psi_{k\tau}=\sum_{l=1}^k\vp_{l\tau}.
$$
The iteration \IterationTauEq\ on $\HO$ can be written as the following system of complex \MA equations
on $\Ho$,
$$
\o_{\psi_{k\tau}}^n=\on e^{f_\o+\frac1\tau\vp_{k\tau}-\mu\psi_{k\tau}}
=\o_{\psi_{(k-1)\tau}}^n e^{(\frac1\tau-\mu)\vp_{k\tau}-\frac1\tau\vp_{(k-1)\tau}}
,\q k\in\NN\eqn\IterationKEEq
$$
(implicit in this equation is also a normalization for $\vp_{k\tau}$ that eliminates the 
ambiguity in passing from an equation on $\HO$ to one on $\Ho$).

We now mention some basic features of the iteration.

The most elementary one is that at each step one gains regularity (two derivatives). This is 
a discrete version of the infinite smoothing property of heat equations.
\TaggPage{SmoothningPropertyPage}

Another distinctive feature of the iteration is that it turns the solution of each type of \MA equation
into the next simplest one. Indeed, to find a \KE metric of negative scalar curvature $-n$ one
needs to solve the equation 
$$
\ovpn=\on e^{\fo+\vp}.\eqn\MANegativeEq
$$
The corresponding time one Ricci iteration requires solving at each step the equation
$$
\ovpn=\on e^{\fo+2\vp}.\eqn\MANegativeTwoEq
$$
Similarly, the Calabi-Yau equation
$$
\ovpn=\on e^{\fo},\eqn\MAZeroEq
$$
is traded for a sequence of equations of the previous type \MANegativeEq, and finally, 
the most difficult equation, 
$$
\ovpn=\on e^{\fo-\vp},\eqn\MAPositiveEq
$$
for a \KE metric of positive scalar curvature $n$, is turned into a sequence 
of Calabi-Yau equations \MAZeroEq\ via the time one iteration, or to a sequence of
equations of the type \MANegativeEq\ for smaller time steps. 

We now discuss the link the iteration creates between classical continuity method paths and the 
heat equation.
In several places in the literature on canonical metrics it is mentioned that the 
continuity method and the heat flow method are morally equivalent 
(see, e.g., \ref{J,Si}).
The following discussion comes to make this statement somewhat more explicit.

Indeed, note that one may also consider the time step of the iteration as a dynamical
parameter and study the continuity method path it defines. More precisely, the time $\tau$ Ricci iteration
is given by a sequence $\{\o_{k\tau}\}_{k\ge1}$ of \K forms in $\HO$ satisfying 
$$
\o_{k\tau}=\o_{(k-1)\tau}+\tau\mu\o_{k\tau}-\tau\Ric\o_{k\tau},\qq \o_0=\o,\eqn\IterationTauTwoEq
$$
for each $k\in\NN$ for which a solution exists.
Now set $k=1$ and consider the path $\{\o_\tau:=\o_{\vp_\tau}\}_{\tau\ge0}$ in $\HO$ (for each $\tau$ for which it exists).
In $\Ho$ we obtain the path
$$
\o_{\vp_{\tau}}^n=\on e^{f_\o+(\frac1\tau-\mu)\vp_{\tau}}
,\q \tau\in(0,\infty).\eqn\ContinuityPathKEEq
$$

Let us compare this path to others that appeared previously in the literature. 

In the case $\mu=1$, when restricted to the segment $\tau\ge 1$ this is just 
a reparametrization of Aubin's path \ref{A2} given by 
$$
\o_{\vp_s}^n=\on e^{f_\o-s\vp_s},\q s\in[0,1], \eqn\AubinPathEq
$$ 
via $s=1-\frac1\tau$. 
Here the solution for \AubinPathEq\  at $s=0$ is given by the Calabi-Yau Theorem. 
Namely, one typically first solves the family of equations introduced by Calabi \ref{C2, (11); Y} 
$$
\o_{\vp_s}^n=\on e^{(s+1)f_\o+c_s},\q s\in[-1,0],\eqn\CalabiPathEq
$$
and then continues to work with the path \AubinPathEq.

For the path \ContinuityPathKEEq\ we still need to invoke the Calabi-Yau Theorem to show
closedness at $\tau\ra 1^-$. 
However, this path can be viewed as a continuity version of the Ricci flow
and has various monotonicity properties that \CalabiPathEq\ does not. When studying the Ricci
iteration for the case $\O=c_1$ this will be useful (note also that this path may be used in place
of \CalabiPathEq\ to prove the Calabi-Yau Theorem).

\FRemm{%
We note that Calabi's path \CalabiPathEq\ can in fact be interpreted as a continuity path arising
from a flow, however not the Ricci flow, see \Ref{OtherDiscreteFlowEq}\
below.
}

In light of this relation to the continuity method, the Ricci iteration is seen to interpolate 
between the continuity method ($\tau=\infty$) and the Ricci flow
($\tau=0$). Cao and Perelman proved that when a \KE metric exists, 
the flow will converge to it in the sense of Cheeger-Gromov. Aubin, Bando-Mabuchi and Yau
proved the analogous result for the continuity method. These results are the main motivation
for \Ref{MainConjecture}.

Next, in the case $\mu=-1$, the continuity path \ContinuityPathKEEq\ 
that arises from the Ricci iteration equation is the same
as the continuity path considered by Tian and Yau in their study of
\KE metrics of negative Ricci curvature on some non-compact manifolds \ref{TY, p. 586}. 
The innovative idea of Tian and Yau was to consider a continuity parameter ``starting from infinity"
observing that along this path one has a uniform lower bound for the Ricci curvature.

Now an open problem concerning the flow equation \RicciFlowEq\ with $\mu>0$ is whether one has a uniform lower bound
for the Ricci curvature depending only on the initial data in the absence of a \KE metric. 
Along the time $\tau$ iteration one does have such a bound, depending
on $\tau$, namely, $\Ric\o_{k\tau}>\frac{\tau\mu-1}\tau\o_{k\tau}$. One possible approach to this problem might
be to show that the Ricci flow stays asymptotically close to the time $\tau$ Ricci iteration for some range of
time steps $\tau$.

\FRemm{
\TaggRemm{DemaillyKollarRemm}%
We remark that when $\mu=1$ another path has been considered previously
by Demailly and Koll\'ar \ref{DK, (6.2.3)}, given by
$\o_{\vp_t}^n=\on e^{tf_\o-t\vp_t},\; t\in[0,1]$. 
As written, this path also does not require to start from a solution to a Calabi-Yau equation.
Yet in order to get openess for it one assumes $\Ric\o>0$ and this involves solving
a Calabi-Yau equation (indeed there is no way to produce a \KE metric without entering $\Hcplus$, and
the Calabi-Yau Theorem amounts to $\Hcplus\ne\emptyset$).
This path can also be explained in terms of a discretization; see \Ref{DemaillyKollarPathEq}.
}

\FRemm{
\TaggRemm{TianZhuRemm}%
Another relation between a continuity path, defined by Tian and Zhu, 
and a discretized flow, this time a modified Ricci flow, will be discussed 
in \Ref{RicciIterationTwistedSection}.
}

\subsection
{Some energy functionals on the space of \K metrics}{} 
\TaggSection{SectionEnergy}%
{\TaggPage{EnergyFunctionalsPage}}%
In this section we will obtain a monotonicity result along the Ricci iteration for a family of energy functionals. 
This result has independent interest, and it seems interesting to compare
it with corresponding studies for the Ricci flow (see \Ref{EkFlowRemm}). 

We briefly recall the pertinent definitions and properties of these energy functionals. For 
more details on this subject we refer to a previous article \ref{R2} and the references therein.

We call a real-valued function $A$ defined on a subset $\Dom(A)$ of 
$\Do\times\Do$ an energy functional if it 
is zero on the diagonal restricted to $\Dom(A)$. By
a Donaldson-type functional, or exact energy functional, we will mean an energy 
functional
that satisfies the cocycle condition $A(\o_1,\o_2)+A(\o_2,\o_3)=A(\o_1,\o_3)$ with each
of the pairs appearing in the formula belonging to $\Dom(A)$
\ref{Do1,M1,T6}. We will occasionally refer to both of these simply as 
functionals and exact functionals, respectively.

The functionals $I,J$, introduced by Aubin \ref{A2}, are defined for each pair 
$(\o,\o_{\vp}:=\o+\i\ddbar\vp)\in\Do\times\Do$ by
$$\eqano{
I(\o,\ovp) & =\Vm\int_M\i\del\vp\w\dbar\vp\w\sum_{l=0}^{n-1}\o^{n-1-l}\w\ovp^{l}
=\Vm\intM\vp(\on-\ovpn),&\eqnno\Ieq\cr
J(\o,\ovp) & =\frac{\Vm}{n+1}\int_M\i\del\vp\w\dbar\vp\w\sum_{l=0}^{n-1}(n-l)\o^{n-l-1}\w\ovp^{l}.
&\eqnno\Jeq
}$$

We say that an exact functional $A$ is bounded from below 
on $U\sseq\Ho$ if for every $\o$ such that $(\o,\ovp)\in\Dom(A)$ and $\ovp\in U$ holds $A(\o,\ovp)\ge C_\o$
with $C_\o$ independent of $\ovp$. 
We say it is proper (in the sense of Tian) on a set $U\sseq\HoG$  if for each $\o\in\HoG$ 
there exists a smooth function 
$\nu_\o:\RR\ra\RR$ 
satisfying $\lim_{s\ra\infty}\nu_\o(s)=\infty$ such that 
$A(\o,\ovp)\ge\nu_\o((I-J)(\o,\ovp))$
for every $\ovp\in U$ \ref{T4, Definition 5.1; T6}. This is well-defined, in other words
depends only on $[\o]$.
Properness of a functional implies it has a lower bound.

Previously \ref{R2} we introduced the following collection of energy functionals for each $k\in\{0,\ldots,n\}$,
$$\eqano{
I_k(\o,\ovp)
= &\
\V\intm\i\del \vp\w\dbar \vp\w\sum_{l=0}^{k-1} \frac{k-l}{k+1}\o^{n-1-l}\w\ovp^l
\cr
= &\ 
\frac\Vm{k+1}\intm\vp(k\on-\sum_{l=1}^k\o^{n-l}\w\ovp^l). &\eqnno\Ikeq
}$$
Note that $I_n=J,\; I_{n-1}=({(n+1)J-I})/n$. 

For each $\mu\in\RR$ define the Ding functional 
$$
F_\mu(\o,\ovp)
=
\cases{ 
\dis-\frac1{n+1}\V\intm\vp\sum_{l=0}^{n}\o^{n-l}\w\ovp^l
-\frac1\mu\log\V\intm e^{f_\o-\mu\vp}\on, 
& for $\mu\ne0$,\cr
\dis-\frac1{n+1}\V\intm\vp\sum_{l=0}^{n}\o^{n-l}\w\ovp^l
+\V\intm \vp e^{f_\o}\on,
& for $\mu=0$.\cr
}\eqn\FFunctionalEq
$$
The critical points of these functionals are \KE metrics \ref{D}. 
However, there is an important difference between the two cases in \FFunctionalEq. While for 
the first the functional is exact, for $\mu=0$ this is not true. This is because 
the second term of $F_0$ is not exact on the space $\Ho$, while the first is.
This rather peculiar phenomenon is reflected also by a property of the generalized Ding functional
(see the end of \Ref{ImprovedMTOEnergyFunctionalsSubS}).

Define the Chen-Tian functionals 
$$
E_k(\o,\ovp)
= \frac{1}{k+1}\V\int_{M\times[0,1]}\big[\mu\o^{k+1}-(\Ric\ovp-\D\dot\vp)^{k+1}\big]\w (\ovp+\dot\vp)^{n-k}\w dt.
\eqn\ChenTianFunctionalsEq
$$
These are well-defined independently of the choice of path and exact \ref{CT}. 
Moreover, in the case $\mu\ne0$, the functionals $E_k$ are related to the
K-energy $E_0$ introduced by Mabuchi \ref{M1} via the following relation:

\FProp{\TaggProp{BMRProp}
Let $\mu\ne0$. 
Let $k\in\{0,\ldots,n\}$. For every $(\o,\ovp)\in \Hmuc\times\Hmuc$,
$$\eqano{
\mu^{k+1}
E_k(\o,\ovp)
& =
\mu
E_0(\o,\ovp)+I_k(\ovp,
\mu
\Ricovp)-I_k(\o,
\mu
\Rico).&\eqnno
}$$
}

\proof
This relation has been previously demonstrated for the case $\mu=1$ \ref{R2, Proposition 2.6}
(recall that the case $\mu=1, k=n$ is a result of Bando and Mabuchi). In general
the same proof goes through by keeping track of the constant $\mu$.\done

Finally, recall the definitions of the following subsets of the space of \K forms on a Fano manifold \ref{R2}:
$$
\eqano{
\calA_k(\o) & =\{\ovp\in\Hc\,:\, E_k(\o,\ovp)\ge0\},&\eqnno\calAkEq
\cr
\calB_k & =\{\ovp\in\Hc\,:\, I_k(\ovp,\Ricovp)\ge0\}.&\eqnno\calBkEq
}$$
When a \KE form $\o$ exists we denote $\calA_k:=\calA_k(\o)$. This is well-defined 
and does not depend on the choice of the \KE form.
We recall that both the sets $\calA_k$ and $\calB_k$ strictly contain the set $\Hcplus$, of \K forms of positive Ricci
curvature (see \Ref{CharacterizationThm} below for a more precise statement).

The following monotonicity result will be useful later. For $E_0, F_\mu$ and $E_1$ with 
$1/\mu=\tau=1$ this was proven before \ref{B,DT,SW}.

\FProp{
\TaggProp{EnergyMonotonicityProp}
(i) The functional $E_0$ is monotonically decreasing along the time $\tau$ iteration 
($\tau>0$) whenever
the initial point is not \KEno. \break
(ii) When $1/\mu=\tau=1$ the same is true for $F_1$, $E_1$, and, 
when the initial metric lies in 
$\calB_k\supset \Hcplus$, also for $E_k, \; k\ge 2$.
}
\proof 
By exactness, it suffices to show monotonicity at each step of the iteration.\newline
(i)
For concreteness % and in view of an application (see \Ref{ImprovedMTOEnergyFunctionalsSubS})
we derive the result only for $\tau=1$ but explicitly compute the the energy decrease
along the iteration (in general see \Ref{TwistedEnergyMonotonicityLemm}, stated for $\mu=1$, that 
works for all $\mu$).
Consider the equation $\o^n_{\vp_1}=\on e^{\fo-(\mu-1)\vp_1}$.
One has \ref{Che1; T4, p. 254; T5, (5.14)}
$$
E_0(\o,\o_{\vp_1})=
\V\intM\log\frac{\o_{\vp_1}^n}{\on}\o_{\vp_1}^n-\mu(I-J)(\o,\o_{\vp_1})+\V\intM \fo(\on-\o_{\vp_1}^n).
$$
First, let $\mu=1$. One has,
$$
E_0(\o,\o_{\vp_1}) = -(I-J)(\o,\o_{\vp_1})+\V\intM\fo\o^n.\eqn\EzeroDecreaseEq
$$
This is nonpositive by the definition of $\fo$ and Jensen's inequality.

When $\mu=-1$ one has 
$$\eqa{
E_0(\o,\o_{\vp_1}) & = \V\intM2\vp_1\o_{\vp_1}^n+(I-J)(\o,\o_{\vp_1})+\V\intM\fo\o^n
\cr & = \V\intM2\vp_1\o_{\vp_1}^n+\V\intM\vp_1(\on-\o_{\vp_1}^n)-J(\o,\o_{\vp_1})+\V\intM\fo\o^n
\cr & = -(I+J)(\o,\o_{\vp_1})+\V\intM(\fo+2\vp_1)\o^n.
}$$
Each term is nonpositive, once again by using the normalization inherent in \MANegativeTwoEq.

When $\mu=0$ one has
$$\eqa{
E_0(\o,\o_{\vp_1}) & = \V\intM\vp_1\o_{\vp_1}^n+\V\intM\fo\o^n
\cr & = -I(\o,\o_{\vp_1})+\V\intM(\fo+\vp_1)\o^n,
}$$
and we may argue as before.

\noindent
(ii)
Using \Ieq-\Jeq\ and \FFunctionalEq\ one has
$$\eqano{
F_1(\o,\o_{\vp_1})
& =-(I-J)(\o,\o_{\vp_1})-\V\intm\vp_1\o_{\vp_1}^n-\log\V\intm e^{\fo-\vp_1}\on&\eqnno\FIterationFirstEq
\cr
& =-(I-J)(\o,\o_{\vp_1})-\V\intm\vp_1\o_{\vp_1}^n-\log\V\intm e^{-\vp_1}\o_{\vp_1}^n.&\eqnno\FIterationSecondEq
}$$
By Jensen's inequality the last two terms combined are nonpositive, and so we conclude by the positivity
of $I-J$. (Alternatively, the iteration stays 
on the submanifold of $\Ho$ defined by the equation 
$\V\intm e^{\fo-\vp}\on=1$ and hence the third term is identically zero, while the second one is negative, again by a special case of Jensen's inequality:
$
1-\V\intM\vp_1\o_1^n\le\V\intM e^{-\vp_1}\o_1^n=\V\intM \o_{2}^n=1$.)

Next, when $1/\mu=\tau=1$ \Ref{BMRProp} gives
$$
E_k(\o,\o_{\vp_1})
=
E_0(\o,\o_{\vp_1})+I_k(\o_{\vp_1},\Ric\o_{\vp_1})-I_k(\o,\Ric\o).\eqn\EkEzeroEq
$$
Next, recall the inequality $I_k\le J$ \ref{R2, (9)}.
Since in our case $\Ric\o_{\vp_1}=\o$ we deduce from \EkEzeroEq\ and \EzeroDecreaseEq\ that
$$
E_k(\o,\o_{\vp_1})\le
-(I-J)(\o,\o_{\vp_1})+J(\o_{\vp_1},\o)-I_k(\o,\Ric\o).
$$
Since $J(\o,\o_{\vp_1})+J(\o_{\vp_1},\o)=I(\o,\o_{\vp_1})$, one has 
$E_k(\o,\o_{\vp_1})\le 0$ if $I_k(\o,\Ric\o)\ge0$, or in other words, if $\o\in\calB_k$.
Finally, note that the subspace $\calB_k$ is preserved under the iteration since, 
in fact, after the first step the iteration will stay in $\Hcplus$.
\done

\FRemm{
An alternative derivation of the second part of (ii) could be to choose a particular path connecting
$\o$ and $\o_1$ and note that each of the contributions has a preferred sign. For example,
choosing Calabi's continuity path \CalabiPathEq\ produces three terms of which 
two are evidently nonnegative. However then
one still needs to manipulate the third term which comes up, 
$-\frac{\Vm}{k+1}\int_M\sum_{i=1}^k {i+1\choose k+1}f_\o(\i\ddbar f_\o)^i\w\o^{n-i}$,
and argue 
that it equals precisely $I_k(\o,\Ric\o)$ and then use the results of \ref{R2} as
above. However to derive (i) for all time steps one needs to use instead the path
\ContinuityPathKEEq\ along which $E_0$ is monotonic, as alluded to after \CalabiPathEq\ above.
}

\FRemm{
\TaggRemm{EkFlowRemm}
Here it is interesting to compare with the Ricci flow. One knows that $F_1,E_0$
are monotonically decreasing along the flow and that as long as
$\Ric\o>-\o$ the same is true for $E_1$ \ref{CT, \S\S3.3, Proposition 4.9}. However an analogous result
is not known along the Ricci flow for $E_k, \; k\ge 2$. In the case that a \KE metric exists one
knows that the flow will converge. One also knows that when restricted to the space $\calB_k$ the functional $E_k$ attains a minimum precisely on the space of \KE metrics, however that outside this space it is not true that 
these functionals are bounded from below on $\Hc$ \ref{R2, \S5} 
(see also \Ref{AnalyticCharacterizationSubS} below). Thus all that is apparent at the 
present moment is that
once the flow stays in $\calB_k$, the functional $E_k$ will eventually decrease,
however even then we do not know whether this will happen monotonically.
}

\subsection
{The Ricci iteration for negative and zero first Chern class}
{\TaggPage{RicciIterationNegativeZeroPage}}% 
\TaggSection{RicciNegativeZeroSection}%
In this section we prove the existence and convergence of the Ricci iteration
in the case that either $c_1<0$ and $\O=-c_1$, or that $c_1=0$ and $\O$
is an arbitrary \K class.

We start with a result that is a simple consequence of the theory
of elliptic complex \MA equations. 
This result is the existence part of \Ref{ConvergenceProp} in the cases under consideration.

\FLem{
Let $(M,\JJJ)$ be a compact \K manifold whose first Chern class is negative or zero. 
When $c_1=0$ denote by $\O$ a \K class; 
otherwise let $\O=\mu c_1$ denote a \K class with $\mu<0$.
Then for any $\o\in\HO$, the time $\tau$ Ricci iteration exists for all $k\in\NN$ and all $\tau\in(0,\infty)$.
}
\proof
It is enough to show existence for one step of the iteration in order
to show the iteration exists for each $k\in\NN$ (by repeating the argument at each step).

Fix $\tau\in(0,\infty)$. 
The existence of $\o_1$ amounts to solving the equation 
$$
\o_{1} =\o_0+\tau\mu\o_{1}-\tau\Ric\o_{1}.
$$
Let $\o_{\vp_1}=\o_1$ with $\vp_1\in\Ho$.
This can be written as a complex \MA equation:
$$
\o_{\vp_1}^n=\on e^{f_\o+(\frac1\tau-\mu)\vp_1},\q \intm \on e^{f_\o+(\frac1\tau-\mu)\vp_1}=V.\eqn\IterationTimetEq
$$
Under the assumption $\mu<1/\tau$, and hence in particular if $\mu\le0$, 
the maximum principle gives an a priori $\Linf$ estimate on $\vp_1$. Then
the work of Aubin and Yau \ref{A1,Y} immediately applies to give higher-order estimates. We conclude
that a unique solution $\o_{\vp_1}\in\HO$ exists. \done

We now turn to the proof of the convergence statement of \Ref{ConvergenceProp} in the cases under consideration.

\bigskip
\noindent
{\it Proof of \Ref{ConvergenceProp} ($\mu\le0$).}
Assume first that $c_1<0$ and let $\O=-c_1$.
We have the following system of \MA equations:
$$
\o_{\psi_{k\tau}}^n
=
\on e^{f_\o+\psi_{k\tau}+\frac1\tau\vp_{k\tau}},\quad k\in \NN.
$$
We first prove an a priori uniform bound, independent of $k$ in an inductive manner.
The first equation reads $\o_{\vp_\tau}^n=\on e^{f_\o+(1+\frac1\tau)\vp_\tau}$. At the maximum of $\vp_\tau$
we have
$\o_{\vp_\tau}\le \o$ and thus $(1+\frac1\tau)\sup\vp_\tau\le -\inf f_\o$. 
A similar argument at the minimum
of $\vp_\tau$ gives $-(1+\frac1\tau)\inf\vp_\tau\le \sup f_\o$.
The second equation reads $\o_{\vp_\tau+\vp_{2\tau}}^n=\o_{\vp_\tau}^n e^{-\frac1\tau\vp_\tau+(1+\frac1\tau)\vp_{2\tau}}$.
The maximum/minimum principle now gives $(1+\frac1\tau)\sup\vp_{2\tau}\le \frac1\tau\sup\vp_\tau$ and 
$-(1+\frac1\tau)\inf\vp_{2\tau}\le -\frac1\tau\inf\vp_\tau$ or 
$\sup\vp_{2\tau}\le -\frac\tau{(1+\tau)^2}\inf f_\o$ and $-\inf\vp_{2\tau}\le\frac\tau{(1+\tau)^2}\sup f_\o$. We then have 
$\sup\psi_{k\tau}\le-\inf f_\o,\; -\inf\psi_{k\tau}\le\sup f_\o $. This uniform bound implies
the existence of an a priori $C^{2,\alpha}$ bound on $\psi_{k\tau}$, independently of $k$.
As a result, by elliptic regularity theory, a subsequence converges to a smooth
solution which we denote by $\psi_\infty$. In fact the convergence is exponentially fast
and there is no need to take a subsequence:
$||\psi_{k\tau}-\psi_{(k-1)\tau}||_{C^{2,\alpha}}\le C \tau (1+\tau)^{-k}$.

Now, by \Ref{EnergyMonotonicityProp}, we notice that unless $\o_0$ is itself
\KEno, the functional $E_0$  
is strictly decreasing along the iteration.
In particular, since  $\o_\infty$ is a fixed point of the iteration it must be
\KEno.

We now consider the case $\mu=0$, for which we have the following
system of equations,
$$
\o_{\psi_{k\tau}}^n
=
\on e^{f_\o+\frac1\tau\vp_{k\tau}},\quad k\in \NN.
$$
We may rewrite this as 
$\o_{\psi_{k\tau}}^n
=
\o_{\psi_{(k-1)\tau}}^n e^{-\frac1\tau\vp_{(k-1)\tau}+\frac1\tau\vp_{k\tau}}$
from which we have $\sup\vp_{k\tau}\le\sup\vp_{(k-1)\tau}\le\ldots\le -\tau\inf f_\o$.
Therefore we have 
$$
||e^{f_\o+\frac1\tau\vp_{k\tau}}||_{L^\infty(M)}\le e^{\oscfoot f_\o},\quad \all\;k\in \NN.
$$
Now, by Yau's work it follows that there exists an a priori $C^{2,\alpha}$ 
bound on $\psi_{k\tau}$, independently of $k$. Combined with the monotonicity result 
it follows, as before, that a subsequence
converges to a \K potential of a \KE metric. 
Moreover, since the \KE metric is unique (in each fixed \K class) \ref{C2} any converging
subsequence will necessarily converge to the same limit point. This then implies
that our original sequence converges to this limit.\done

\subsection
{The Ricci iteration for positive first Chern class}
\TaggSection{RicciIterationPositiveSection}%
{\TaggPage{RicciIterationPositivePage}}%
We turn to the study of the iteration on Fano manifolds that, as noted in the Introduction, is
our main motivation for introducing the Ricci iteration. Most of the applications described
in \Ref{ApplicationsSection} are for this class of manifolds.

We first 
introduce an operator that arises very naturally although it seems to 
have not been defined previously in the literature. 
It exists and is well-defined by the Calabi-Yau Theorem \ref{Y}.

\FDef{
\TaggDef{InverseRicciDef}
Define the
inverse Ricci operator $\Ric^{-1}:\Dc\ra\Hc$ by letting 
$\Ric^{-1}\o:=\o_\vp$ with $\ovp$ the unique \K form in $\Hc$ 
satisfying $\Ric\ovp=\o$. Similary denote higher order iterates of this operator
by $\Ric^{-l}$ for each $l\in\NN$. Let $\Ric^0:=\h{\rm Id}$ denote the identity operator.
}
There exists a generalization of this operator to any \K manifold 
(\Ref{GeneralizedInverseRicciDef}). For another direction in which this
operator may be generalized see \Ref{ATwistedInverseRicciDef}.

We then see that the dynamical system corresponding to the 
time one Ricci iteration on a Fano manifold with $\mu=1$
is nothing but the evolution of iterates of the inverse Ricci operator, 
$$
\o_l=\Ric^{-l}\o_0.
$$

The following result concerns the ``allowed" time steps in the iteration for any Fano manifold
and is well-known.
Note that unlike in the previous, unobstructed, cases, the allowed range for the time step
is restricted unless an analytic ``semi-stability" condition holds.

Define 
$$
\tau_M(G)=\sup\{\;\,t \,:\, \h{\ \ContinuityPathKEEq\ has a solution for each } \tau\in(0,t)
\h{\ and } \o\in\HcG  \}. \eqn\TauInvariantEq
$$
By definition this is a holomorphic invariant. Recall also the definition of Tian's invariants
\ref{T1,T3}
$$\eqano{
\alpha_M(G)& =\sup\{\; a \,:\, \sup_{\vp\in\calfootH_\o(G)}\V\intm e^{-a(\vp-\sup\vp)}\on<\infty  \}, &\eqnno\TianInvariantEq\cr
\beta_M(G) & =\sup\{\; b \,:\, \Ric\o\ge b\o, \;\o\in\HO(G)  \}, &\eqnno\TianInvariantTwoEq\cr
}$$
where in \TianInvariantEq\ $\o$ is any element of $\HcG$.

\FLem{
\TaggLemm{FanoIterationTimeStepsLemm}
Let $(M,\JJJ)$ be a Fano manifold and let $G$ be a compact subgroup of $\AutMJ$. \newline
(i) For any $\o\in\HcG$, the time $\tau$ Ricci iteration exists for all 
$k\in\NN$ and all $\tau\in[0,\tau_M(G))$. One has 
$\frac1{1-\beta_M(G)}\ge\tau_M(G)\ge \big|\frac 1{\max\{1-(n+1)\alpha_M(G)/n,0\}}\big|>1$.\newline
(ii) Assume that $E_0$ is bounded from below on $\HcGplus$. 
Then $\tau_M(G)=\infty$.
}

\proof (i) By the Calabi-Yau Theorem $\tau_M(G)\ge1$.
According to Tian \ref{T1} the path \AubinPathEq\ exists for each 
$s\in[0,(n+1)\alpha_M(G)/n)\cap[0,1]$
whenever $\o\in\Hc(G)$. Note that $\tau=1/(1-s)$ and that $\alpha_M(G)>0$.
\newline
(ii) This is equivalent to a result of Bando and Mabuchi \ref{BM, Theorem 5.7}.
\done

Combined with \Ref{BMDTThm} we therefore obtain the existence part of 
\Ref{ConvergenceProp} for $\mu>0$.

\FCor{
Let $(M,\JJJ)$ be a \KE Fano manifold.
Then for any $\o\in\Hc$, the time $\tau$ Ricci iteration exists for all $k\in\NN$ and all $\tau\in(0,\infty)$.
}

We now turn to the proof of the remaining part of \Ref{ConvergenceProp}. The argument is essentially due
to Bando and Mabuchi \ref{BM}.
We assume for simplicity, as in the statement of the theorem, that $\autMJ=\{0\}$ (for the additional details necessary for the 
general case we refer to \ref{R5} (see also \Ref{RicciIterationTwistedS})).
For the case $\mu=1$ we are solving the system of equations
$$
\o_{\psi_{k\tau}}^n
=
\on e^{f_\o-\psi_{k\tau}+\frac1\tau\vp_{k\tau}},\quad k\in \NN.\eqn\RicciIterationPositiveEq
$$
Let $G_{k\tau}$ be a Green function for $-\D_{k\tau}=-\D_{\dbar,\o_{\vp_{k\tau}}}$ 
satisfying $\intM G_{k\tau}(\cdot,y)\o_{\psi_{k\tau}}^n(y)=0$. Set $A_{k\tau}=-\inf_{M\times M} G_{k\tau}$.
Since $-n<\D_0\psi_{k\tau}$ and $n>\D_{k\tau}\psi_{k\tau}$ the Green formula gives
$$
\eqano{
\psi_{k\tau}(x)-\VintM \psi_{k\tau}\o_0^n = & -\VintM G_0(x,y)\D_0\psi_{k\tau}(y)\o_0^n(y)\le nA_0,\cr
\psi_{k\tau}(x)-\VintM \psi_{k\tau}\o_{\psi_{k\tau}}^n= & -\VintM G_{k\tau}(x,y)\D_{k\tau}\psi_{k\tau}(y)\o_{\psi_{k\tau}}^n(y)\ge -nA_{k\tau}.\cr
}
$$
Hence
$$
\osc\psi_{k\tau}\le n(A_0+A_{k\tau})+I(\o_0,\o_{\psi_{k\tau}}).\eqn\OscEq
$$

Since by \Ref{BMDTThm} (ii) $E_0$ is proper on $\Hc$ in the sense of Tian, 
if $E_0(\o,\cdot)$ is uniformly bounded from above on a subset of $\Hc$ so is $I(\o,\cdot)$. 
By the monotonicity of $E_0$ along the iteration  
we conclude that $I(\o,\o_{\psi_{k\tau}})$ is uniformly bounded independently of $k$. 

It remains to bound $A_{k\tau}$. This can be done using a 
special case of Bando and Mabuchi's Green's function estimate that we now state.
\FThm{
\TaggThm{BMGreenFunctionThm}
\ref{BM, Theorem 3.2}
Let $(N,h)$ be a connected compact closed Riemannian manifold of nonnegative 
Ricci curvature. Let $G_h$ denote the Green function of $d^{\star_h}\circ d+d\circ d^{\star_h}$
satisfying $\int_N G(x,y)dV_h(y)=0$ for each $x\in N$ and let 
$A_h=-\inf_{M\times M} G_h$. Then
$$
A_h\le c_n \frac {\diam(N,h)^2}{\vol(N,h)},
$$
with $c_n$ depending only on $n$.
}

Now, along the iteration it holds $\Ric\o_k>(\tau-1)\o_k>0$. 
By Myers' Theorem \ref{P, p. 245} then 
$$
\diam(M,\o_k)^2\le\pi^2(2n-1)/(\tau-1).\eqn\MyersDiameterBoundEq
$$ 
Set $\tilde\psi_{k\tau}:=\psi_{k\tau}-\V\intm\psi_{k\tau}\on$.
Now \Ref{BMGreenFunctionThm} and equations \OscEq-\MyersDiameterBoundEq\ yield the 
estimate $||\tilde\psi_{k\tau}||_{\Linf}\le C$. As in \Ref{RicciNegativeZeroSection}, the general theory 
of \MA equations now provides for uniform higher derivative estimates. 
We may thus extract a converging subsequence from $\{\tilde\psi_{k\tau}\}_{k\ge0}$. 
Thanks to the monotonicity of $E_0$
it must converge to a \K potential for a \KE metric. Since such a metric
is unique \ref{BM, Remark 9.3} the same argument as before gives the convergence of the full
orbit of the Ricci iteration.\done

\subsection
{The \KR flow and the Ricci iteration for a general \K class}
\TaggSection{RicciIterationGeneralSection}%
{\TaggPage{RicciIterationGeneralPage}}%
A natural question is whether on an arbitrary \K manifold 
one may define an iteration scheme generalizing the Ricci iteration. To answer this question of course
one first needs to generalize the Ricci flow itself.
In this section we recall one such possibility. 
We end with a conjecture regarding the convergence of this iteration.

A flow on the space of \K forms $\HO$ can be considered as an integral
curve of a vector field on this space. A vector field $\chi$ on $\HO$ is an assignment 
$\o\mapsto \chi_\o\in \CinfM/\RR$. 
The Ricci flow describes the dynamics of minus the Ricci potential
vector field $-f$. Recall that
the vector field $f$ is the assignment $\o\mapsto f_\o$ with $\fo$ defined by 
$\Ric\o-\mu\o=\i\ddbar\fo$, $\,\mu\in\RR$, where $\mu \O=c_1$. 

The Ricci iteration in turn can be thought of
as a piecewise linear trajectory in $\HO$ induced from the Ricci potential vector field $-f$ 
and approximating its integral curves. 

Motivated by this one is naturally led to extend the definition of the Ricci flow \RicciFlowEq\ to an arbitrary \K
manifold, simply by defining the flow lines to be integral curves of minus the Ricci potential vector field $-f$
on $\HO$, with $\O$ an arbitrary \K class. 
Recall that the Ricci potential is defined in general by 
$\Ric\o-H_\o\Ric\o=\i\ddbar\fo$. The resulting flow equation can be written as
$$
\eqa{
\frac{\pa \o(t)}{\pa t} & =-\Ric\o(t)+H_t\Ric\o(t),\quad t\in\RR_+,\cr
\o(0) & =\o,
}\eqn\KRFEq$$
for each $t$ for which a solution exists in $\HO$
(throughout subscripts are meant to indicate that the relevant object corresponds to the metric 
indexed by that subscript). This flow, introduced by Guan, 
is part of the folklore in the field although it has not been much studied.\note{%
\cmreight It seems that Guan first considered this flow in unpublished work in
the 90's (see references to\break%
\vglue-0.5cm 
\hglue-\parindent\cmreight 
\reffoot{G1}). After completing this article I also became aware, thanks to
G. Sz\'ekelyhidi, of a recent preprint\break %  
\vglue-0.5cm 
\hglue-\parindent\cmreight \reffoot{G2}  posted by Guan on his webpage in which
this flow is studied. We hope that the elementary\break %  
\vglue-0.5cm 
\hglue-\parindent\cmreight 
 discussion in this section is still
of some interest even though it was written before learning of \reffoot{G1,G2}.\break %  
\vglue-0.5cm 
\hglue-\parindent\cmreight  
For a different but related flow see \reffoot{S}.
}

Corresponding to this flow we introduce the following dynamical system on $\HO$ which 
generalizes \Ref{RicciIterationDef}.

\FDef
{
\TaggDef{RicciIterationGeneralizedDef}%
Given a \K form $\o\in\HO$ let the time $\tau$ Ricci iteration
be 
the sequence of forms $\{\o_{k\tau}\}_{k\ge0}$,
satisfying the equations
$$
\eqa{
\o_{k\tau} & =\o_{(k-1)\tau}+\tau H_{k\tau}\Ric\o_{k\tau}-\tau \Ric\o_{k\tau},\quad   k\in\NN,\cr
\o_0 & =\o,
}\eqn\RIterationEq$$
for each $k\in\NN$ for which a solution exists in $\HO$. 
}

\noindent
As in \Ref{SectionRicciIteration}, setting $k=1$ and varying $\tau$ defines a continuity path
that is of independent interest.

An observation that goes back to Calabi characterizes the equilibrium state of the flow 
and the iteration.
\FLem{\ref{C1, Theorem 1}
The Ricci form of a \K metric is a harmonic representative of $c_1$ with respect to the metric
if and only if its scalar curvature is constant.
}
\proof
One has
$$
n\Ric\o\w\o^{n-1}=\tr_\o\Ric\o\;\on=s(\o)\on.
$$
Since $\o$ is a harmonic representative of its class, we see that $s(\o)$ is harmonic, i.e., constant,
if and only if $\Ric\o$ is. \done

An infinitesimal automorphism $X\in\autMJ$ naturally induces a vector field $\psi^X$ on $\HO$
given by 
$$
\psi^X:\o\mapsto\psi^X_\o\in\CinfM/\RR, \q \h{where $\L_X\o=\i\ddbar\psi^X_\o$.}\eqn\PsiXVectorFieldEq
$$ 
Recall the following generalization of the notion of a constant scalar curvature \K metric, due to Guan. Alternatively it may
be seen as a generalization of the notion of a \KR soliton to an arbitrary class.
\FDef{
\TaggDef{KRSolitonDef}\ref{G1}
Let $X\in\autMJ$. A \K metric $\o$ will be called a \KR soliton if it satisfies 
$$
\Ric\o-H_\o\Ric\o=\Lie_X\o.\eqn\KRSolitonEq
$$
Equivalently, if the vector field $\psi^X-f$ on $\H_{[\o]}$ has a zero at $\o$. 
}

Motivated by the results for \KE manifolds we believe the following conjecture should hold.

\FConj{
Let $(M,\JJJ)$ be a compact closed \K manifold, and assume that there
exists a constant scalar curvature \K metric representing the class $\O$.
Then for any $\o\in\HO$, the \KR flow \KRFEq\ and the
Ricci iteration \RIterationEq\ exist 
and converge in an appropriate sense to a constant scalar curvature metric. 
}

Similarly, we believe an analogous result should hold for \KR solitons \KRSolitonEq\ using
the twisted constructions of \Ref{RicciIterationTwistedSection}.

\subsection
{Another flow and the inverse Ricci operator for a general \K class}
{\TaggPage{AnotherFlowPage}}%
\TaggSection{AnotherFlowSection}%
Our purpose in this section is to explain why the inverse Ricci operator---that 
appeared as a very singular iterative construction for anticanonically polarized
Fano manifolds---is in fact a special case of a more general construction on any \K manifold. 
This gives another application of our approach explained in the Introduction since it involves a
discretization of another geometric flow equation.

To that end, given a \K form $\o$ let us consider the flow equations
$$
\eqa{
\frac{\pa \Ric\o(t)}{\pa t} & =-\Ric\o(t)+H_t\Ric\o(t),\quad t\in\RR_+,\cr
\o(0) & =\o,
}\eqn\AnotherFlowEq
\TaggEq{AnotherFlowEq}
$$
for each $t$ for which a solution exists.

The following brief and informal discussion comes to motivate this definition.  
Consider the case when the first Chern class is definite 
($\mu\in\RR\setminus\{0\}$ with $\O=\mu c_1$), or zero ($\O$ is arbitrary),
and take $\o\in\HO$. The evolution equation then becomes
$$
\frac{\pa \Ric\o(t)}{\pa t}  =-\Ric\o(t)+\mu\o(t).\eqn\AnotherFlowDefiniteCaseEq
$$
Assume momentarily that the flow preserves the \K class and that it exists on some time interval $[0,T]$. 
Then on the level of potentials it can be written as
$$
-\D_t\dot\vpt=\log\frac{\ovptn}{\on}+\mu\vpt-\fo+a_t,\q \vp_0=\h{const}, \eqn\AnotherFlowPotentialEq
$$
or as a \MA equation
$$
\ovptn=\on e^{\fo-\mu\vpt-\D_t\dot\vpt-a_t},\eqn\AnotherFlowPotentialMAEq
$$
with $a_t$ a certain normalizing constant.
Set $u:=\D_t\dot\vpt$. 
A time derivative of \AnotherFlowPotentialEq\ gives
$$
\frac{du}{dt}=-u+\mu G_t u+b_t,
$$
with $b_t$ another normalizing constant.
One may show that
$$
||u||_{\Linf(M\times[0,T])}<Ce^{-t},
$$
when $\mu\le0$ and that $||u||_{\Linf(M\times[0,T])}<Ce^{(\mu/\lambda_1(t)-1)t}$,
when $\mu>0$, where $\lambda_1(t)$ is the first nonzero eigenvalue of $-\D_t$. The constant
$C$ depends a priori on $t$. 
Going back to \AnotherFlowPotentialMAEq\ one may show an a priori
estimate $||\vpt||_{\Linf(M\times[0,T])}<C_1$, with $C_1$ depending
only on $\o$, whenever $\mu\le 0$. This then implies a priori estimates on
higher order derivatives. Finally, take a converging subsequence. Along this subsequence
$\lambda_1$ is uniformly bounded away from zero. Going back 
to the exponential decay of $u$ we apply uniform Schauder estimates to conclude that $\dot\vpt$
is uniformly decaying. It then follows that the limit is a \KE metric. By uniqueness of the metric
one then argues, as earlier on, that the flow itself converges exponentially fast to
a \KE metric. On the other hand, the case $\mu>0$ would require more work, quite likely 
in the spirit of the corresponding result for the Ricci flow \ref{Ca,TZ4} (cf. also \ref{A2,BM}).

Motivated by this discussion, we introduce the following dynamical system on $\HO$ obtained as 
the time one Euler method for this flow:
$$
\eqa{
\Ric\o_{k+1} & =H_{k}\Ric\o_{k},\quad   k\in\NN,\cr
\o_0 & =\o.
}\eqn\OtherIterationEq
\TaggEq{OtherIterationEq}
$$
It can be thought of as describing the dynamics of a generalized inverse Ricci operator.
This motivates the following definition, generalizing \Ref{InverseRicciDef} to an 
arbitrary \K manifold.

\FDef{
\TaggDef{GeneralizedInverseRicciDef}
Define the
inverse Ricci operator $\Ric_\O^{-1}:\HO\ra\HO$ by letting 
$\Ric_\O^{-1}\o:=\o_\vp$ with $\ovp$ the unique \K form in $\HO$ 
satisfying $\Ric\ovp=H_\o\Ric\o$. Similary we denote higher order iterates of this operator
by $\Ric_\O^{-l}$ for each $l\in\NN$.
}

Calabi-Yau manifolds are singled-out as those manifolds for which this operator is a constant
map. In general the dynamics of this operator seem intriguing.

We end this section with two remarks regarding continuity method paths induced from the 
flow \AnotherFlowDefiniteCaseEq, directly continuing the discussion in \Ref{SectionRicciIteration}.
First, it is interesting to note that discretizing this flow for time steps $\tau\in[0,1]$ 
gives rise to the well-known continuity path of the Calabi-Yau Theorem (here $\mu=0$)
introduced by Calabi \ref{C2, (11)},%
\note{\cmreight To obtain this path in the equivalent setting of the search for a \K metric 
with prescribed\break 
\vglue-0.5cm
\hglue-\parindent\cmreight Ricci
form, one considers the flow obtained by replacing the harmonic projection term 
in \AnotherFlowEq\ by a\break 
\vglue-0.5cm
\hglue-\parindent\cmreight  prescribed form representing $\meight c_1$.}
$$
\Ric\o_{\vp_\tau}-\Ric\o=-\tau\Ric\o\q \Longrightarrow\q e^{\tau f_\o+d_\tau}\on=\o_{\vp_\tau}^n,
\eqn\OtherDiscreteFlowEq
\TaggEq{OtherDiscreteFlowEq}
$$
with $d_\tau=-\log\V\intm e^{\tau f_\o}\on,\; \tau\in[0,1]$.

The K-energy decreases along this path, however not monotonically in general.
This is in contrast to the continuity path arising from the Ricci iteration and 
fits-in well with what we would expect: 
the former arises from the Euler method (as opposed to the backwards Euler method) 
and so one does not expect monotonicity, nor convergence for large enough time steps.

Also, we remark that the backwards Euler method of the same evolution equation
\AnotherFlowDefiniteCaseEq\ yields the continuity path
$$
\o_{\vp_\tau}^n=\on e^{\frac \tau{1+\tau}(f_\o-\mu\vp)}, \q \tau\ge 0,\eqn
$$
\TaggEq{DemaillyKollarPathEq}%
that coincides in the case $\mu=1$, after reparametrization, with the continuity path used by Demailly and Koll\'ar alluded to earlier (\Ref{DemaillyKollarRemm}).

\subsection
{The twisted Ricci iteration and a twisted inverse Ricci operator}
\TaggSection{RicciIterationTwistedSection}%
\TaggS{RicciIterationTwistedS}%
\TaggPage{RicciIterationTwistedPage}%
When searching
for canonical metrics,
the presence of continuous symmetries has traditionally required additional analysis. Although
the arguments are very similar to the previous sections, there are certain differences.
In this section we merely introduce some of the dynamical constructions relevant to this case which
will be further used and studied in the sequel \ref{R5} in the setting of convergence towards
\KE metrics with continuous symmetries and \KR solitons (or multiplier Hermitian structures). We also state a monotonicity result that will
be used in \Ref{NadelQuestionSubS}.

In the presence of holomorphic vector fields one oftentimes modifies the flow equation by a time-dependent
family of automorphisms \ref{CT,TZ4}. More generally, one may study the dynamics of a perturbation of the
vector field $-f$ by an arbitrary vector field $\chi$. Adapting the point of view of either 
\Ref{RicciIterationGeneralSection} or \Ref{AnotherFlowSection} yields two ways to obtain
discrete dynamics.
The following definition corresponds to the former.

\FDef{
\TaggDef{GeneralizedTwistedRicciIterationDef}
Given a vector field $\chi:\o\mapsto\chi_\o\in\CinfM/\RR$ 
on $\HO$ define the $\chi$-twisted time $\tau$ Ricci iteration
to be the sequence of forms $\{\o_{k\tau}\}_{k\ge0}$ satisfying the equations
$$
\eqa{
\o_{k\tau} & =\o_{(k-1)\tau}+\tau H_{k\tau}\Ric\o_{k\tau}-\tau \Ric\o_{k\tau}+\tau\i\ddbar\chi_{\o_{k\tau}},\quad   k\in\NN,\cr
\o_0 & =\o,
}\eqn\TwistedRIterationEq$$
for each $k\in\NN$ for which a solution exists in $\HO$. 
}

The construction in \Ref{RicciIterationGeneralSection} corresponds
to the zero vector field. % while the vector field $f$ yields the identity operator. 
The case $\chi=\psi^X$, with $X$ an infinitesimal automorphism, 
will be useful when studying convergence towards solitons.

When $\O=c_1, \tau=1$ this iteration takes on a special form, giving a certain
generalized inverse Ricci operator (cf. \Ref{InverseRicciDef}).

\FDef{
\TaggDef{ATwistedInverseRicciDef}
Given a vector field $\chi:\o\mapsto\chi_\o\in\CinfM/\RR$ 
on $\Hc$ 
define the $\chi$-twisted Ricci
operator $\Ric_{\chi}:\Hc\ra\Dc$ by letting $\Ric_{\chi}\o:=\Ric\o-\i\ddbar\chi_\o$. 
Define the $\chi$-twisted inverse Ricci operator
$\Ric_{\chi}^{-1}:\Hc\ra\Hc$ by letting
$\Ric_{\chi}^{-1}\o:=\o_\vp$ whenever there exists a unique \K form $\ovp$ in $\Hc$ 
satisfying $\Ric_\chi\ovp=\o$. 
Denote higher-order iterates of these operators by $\Ric_\chi^l$ for $l\in\ZZ$, setting
$\Ric_\chi^0:=\h{\rm Id}$.
}

Recall that the Bakry-\'Emery Ricci form associated to a pair $(\o,a)\in\HO\times\CinfM$
is the form $\Ric\o-\i\ddbar a$, that can viewed as the Ricci form of the \K manifold 
$(M,\JJJ)$ equipped with a \K form whose top exterior product equals
$e^{2\pi a}\on$ \ref{BE, (4b)}. The twisted Ricci operator is thus an assignment of a Bakry-\'Emery Ricci form to
each \K form determined by a vector field on $\HO$. The simplest examples include the zero vector field and the Ricci 
potential vector field that yield the Ricci operator and the identity operator, respectively.
Note that the fixed points of %%$\Ric_{\O,\chi}^{-1}$ 
the twisted Ricci operator are for certain choices of $\chi$
the multiplier Hermitian structures defined by Mabuchi \ref{M3}.  
The twisted inverse Ricci operator is not defined for general $\chi$, however it is for some geometrically
significant vector fields.
Assume that $X$ belongs to a reductive Lie subalgebra of $\autMJ$ and that the 
one-parameter subgroup $T_{\h{\notsosmall J} X}$ generated by $\JJJ X$ is 
a compact torus in $\AutMJ$.
When $\chi=\psi^X$ the operator $\Ric^{-1}_{\psi^X}$ restricted to $\Hc(T_{\h{\notsosmall J} X})$ exists and 
is well-defined according to a  theorem of Zhu \ref{Z}. More generally, this is still true when $\chi$ is a smooth function of
$\psi^X$ under some assumptions \ref{M3}.

First, continuing the discussion of \Ref{SectionRicciIteration} (\Ref{TianZhuRemm})
observe that when $\O=c_1$ and $\chi=\psi^X$ the continuity method
path obtained by setting $k=1$ and letting $\tau$ vary in the segment $[1,\infty)$ coincides with the Tian-Zhu
continuity path \ref{TZ2, (1.4)}
$$
\o_{\vp_s}^n=\on e^{f_\o-\psi^X_{\o_{\vp_s}}-s\vp_s},\q s\in[0,1], \eqn\TianZhuPathEq
$$
via the reparametrization $s=1-\frac1\tau$, discretizing the $\psi^X$-twisted \KR flow \ref{TZ4, (4.4)}
$$
\eqa{
\frac{\pa \o(t)}{\pa t} & =-\Ric\o(t)+\o(t)+\L_X\o(t),\quad t\in\RR_+,\cr
\o(0) & =\o\in\Hc(T_{\h{\notsosmall J} X}).
}\eqn\ModifiedRicciFlowEq
$$
In fact, more generally Mabuchi's continuity path \ref{M3, (5.1.4)} in the context of multiplier Hermitian structures
is obtained in the same manner from \TwistedRIterationEq\ as a result of discretizing the corresponding twisted
\KR flow.

We now discuss briefly the special case of
\KR solitons. This is mainly done for the sake of concreteness since, due to the work of Mabuchi, the relevant
computations go through also for general multiplier Hermitian structures.

In their study of \KR solitons on Fano manifolds Tian and Zhu introduced a twisted 
version of the functional %% $F_1$ and 
$E_0$ \ref{TZ3}. To define it we first recall some relevant 
facts \ref{F1, \S2.4; M3,TZ2}.
Given $X\in\autMJ$, let $L^{\psi^X}_\o$ denote the elliptic operator
$L^{\psi^X}_\o\phi:=\D_\o\phi+X\phi$. This operator is self-adjoint with respect to the 
$L^2(M,e^{\psi^X_\o}\on)$ inner product denoted by $\langle\;\cdot\;,\;\cdot\;\rangle_{\psi^X}$.
The vector field $\psi^X$ on $\Hc$ induces a vector field on the space of \K potentials (that
we still denote by the same notation) by decreeing that $\V\intm e^{\psi_\o^X}\on=1$ for each
$\o\in\Hc$. One then has $\psi^X_{\ovp}=\psi^X_\o+X\vp$ since 
$\ddt\V\intm e^{\psi_{\o}^X+X(t\nu)}\o_{t\nu}^n=\langle L_{\o_{t\nu}}^{\psi^X}\nu,1\rangle_{\psi^X}=0$. Define
a functional on $\Hc(T_{\h{\notsosmall J} X})\times\Hc(T_{\h{\notsosmall J} X})$ by

$$
\eqano{
E_0^{\psi^X}(\o,\ovp) & 
= \V\int_{[0,1]}\langle\dot\vpt,L_{\ovpt}^X(\psi^X_{\ovpt}-\fovpt)\rangle_{\psi^X} dt
&\eqnno\ModifiedEzeroFunctionalEq
}
$$
This functional is well-defined independently of a choice of path and exact. Its critical points
are \KR solitons.

\FLem{
\TaggLemm{TwistedEnergyMonotonicityLemm}%
Assume that $X$ belongs to a reductive Lie subalgebra of $\autMJ$ and that the 
one-parameter subgroup $T_{\h{\notsosmall J} X}$ generated by $\JJJ X$ is 
a compact torus in $\AutMJ$.
The functional 
 $E_0^{\psi^X}$  
 is monotonically decreasing along the 
$\psi^X$-twisted time $\tau$ Ricci iteration for each $\tau>0$ for which the iteration exists.}
\proof 
Let $\o\in\Hc(T_{\h{\notsosmall J} X})$.
Whenever the $\psi^X$-twisted time $\tau$ Ricci iteration exists the same is
true for smaller time step iterations. Hence the 
continuity path 
$$
\ovptn=\on e^{\fo-\psi^X_{\ovpt}+(\frac1t-1)\vpt},\q t\in[0,\tau],\eqn\TwistedCMPathEq
$$ 
exists.
Differentiating equation \TwistedCMPathEq\ gives
$(L_{\ovpt}^{\psi^X}+1-\frac1t)\dot\vpt=-\frac1{t^2}\vpt$. Hence
one has
$$\eqa{
E_0^{\psi^X}(\o_0,\o_\tau) & =
\V\int_{[0,\tau]}\frac1t\langle\dot\vpt,L_{\ovpt}^{\psi^X}\vpt\rangle_{\psi^X} dt\cr
& = -\V\int_{[0,\tau]}t\langle L_{\ovpt}^{\psi^X}\dot\vpt,(L_{\ovpt}^{\psi^X}+1-\frac1t)\dot\vpt\rangle_{\psi^X} dt\le0.
}$$
When $\tau\le 1$ the last inequality is a just a consequence of the ellipticity of $L^{\psi^X}_{(\,\cdot\,)}$. 
When $\tau>1$ it follows since
$L_{\ovpt}^{\psi^X}+1-\frac1t$ is still elliptic \ref{TZ2, Lemma 2.2 (ii)}.
\done

\subsection
{Some applications}
\TaggSection{ApplicationsSection}
{\TaggPage{ApplicationsPage}}%
In this section we describe several applications of the Ricci iteration and the
inverse Ricci operator
to some classical objects and problems in \K and conformal geometry. 

\subsubsection{The Moser-Trudinger-Onofri inequality on 
the Riemann sphere and its higher dimensional analogues}
\TaggSubsection{MTOSubsection}%
\TaggSubS{MTOSubS}%
We recall some notions from \ref{R2} and explain how the results there on the Moser-Trudinger-Onofri
inequality can be rephrased in terms
of the inverse Ricci operator.
This sheds new light on our discussion there and at the same time expands it (this was omitted
from \ref{R2} for the sake of brevity).

{\TaggPage{MTOPage}}

Let $\oFSc$ denote the Fubini-Study form of constant Ricci curvature $c$ on $(S^2,\JJJ)$, 
the Riemann sphere, given locally by
$$
\oFSc=\frac\i{c\pi}\frac{dz\w\dbz}{(1+|z|^2)^2}.
$$
Here $V=\int_{S^2}\oFSc=c_1([M])/c=2/c$. For
$c=1/2\pi$ it is induced from restricting the Euclidean metric on $\RR^3$ to 
the radius $1$ sphere. Denote by $W^{1,2}(S^2)$ the space of functions on $S^2$
that are square-summable and so is their gradient (with respect to some Riemannian 
metric).
The Moser-Trudinger-Onofri inequality states:
\FThm{\ref{Mo,O,Tr}
\TaggThm{MTOThm}
 For
$\o=
\o_{\hbox{\small FS},2/V}
$ and 
any function $\vp$ on $S^2$ in $W^{1,2}(S^2)$ one has
$$
\V\int_{S^2} e^{-\vp+\V\int_{S^2}\vp\o}\o\le e^{\V\int_{S^2}\ha\i\del\vp\w\dbar\vp}.\eqn\MTClassicalineq
$$
Equality holds if and only if $\ovp$ is the pull-back of $\o$ by a M\"obius transformation.
}
\noindent

Several proofs of this classical result have been given in the literature, and we list here
the ones we are aware of, chronologically: Onofri \ref{O}, Hong \ref{H}, Osgood-Phillips-Sarnak \ref{OPS},
Beckner \ref{Be}, Carlen and Loss \ref{CL1,CL2}, Ghigi \ref{Gh} (for more background we refer to Chang \ref{Ch}).
All of these proofs use crucially some symmetrization/rearrangement arguments that reduce the problem
to a single dimension.
Previously we gave a new proof of this inequality coming from \K geometry \ref{R2}. 
At the same time we also formulated an optimal (in a sense to be clarified below) extension of it
to higher-dimensional \KE manifolds of positive scalar curvature, extending the work of Ding and Tian. 

A function satisfies \MTClassicalineq\ if and only if
$$
F_1(\o,\ovp)=\V\int_{S^2}\ha\i\del\vp\w\dbar\vp-\V\int_{S^2}\vp\o-\log\V\int_{S^2} e^{-\vp}\o\ge0.
$$
This functional was studied already by Berger and Moser \ref{Ber,Mo}. Moser, extending work
of Trudinger, was able to show that $F_1(\o,\ovp)\ge -C$. Then, Onofri showed that $C=0$ and characterized
the cases of equality.

Aubin first suggested a connection between the classical inequality
\MTClassicalineq\ and \KE manifolds \ref{A2, (4)}. Following this, Ding \ref{D} showed
how to generalize the functional $F_1$ to higher-dimensional Fano 
manifolds---see Equation \FFunctionalEq---using Aubin's functional $J$.
Using this observation, and
modifying the proof of a fundamental result of Bando and Mabuchi 
concerning the boundedness of the K-energy,
Ding and Tian proved \Ref{MTOThm} for those functions that
belong to the subspace 
$\calH_\o\subset W^{1,2}(S^2)$. 
We state both results and their corollary. The corollary
is Ding and Tian's restricted generalization\note{\cmreight This term 
is meant to emphasize that this generalized a restricted version of the classical\break 
\vglue-0.5cm
\hglue-\parindent\cmreight inequality.}
of the Moser-Trudinger-Onofri inequality to higher-dimensional
\KE manifolds.

\FThm{\TaggThm{BMBound}\ref{BM, Theorem A; B, Theorem 1; DT, Theorem 1.1} 
Let $(M,\JJJ,\o)$ be a \KE Fano manifold. Then
$E_0(\o,\ovp), F_1(\o,\ovp) \ge0$ for all $\ovp\in\Hc$
with equality if and only
if $\ovp=h^\star\o$ with $h\in\Aut(M,\JJJ)_0$.
}

\vglue-0.3cm

\FCor{
\TaggCor{DTBoundednessThm}
\unskip\unskip\unskip\ref{DT}
Let $(M,\JJJ,\o)$ be a \KE Fano manifold with $\o\in\Hc$.
Then for each $\vp\in\H_\o$ holds
$$
\V\intm e^{-\vp+\V\intm\vp\on}\on
\le 
e^{J(\o,\ovp)}.\eqn\GeneralizedMTOeq
$$
Equality holds if and only if $\ovp$ is the pull-back of $\o$ by a holomorphic transformation.
}

One should note that the subspace of \K potentials can be considered
as a rather small ``ball" sitting inside $\Cinf(S^2)\subset W^{1,2}(S^2)$ since in general a large enough multiple
of an element of $\calH_\o$ will no longer belong to $\calH_\o$. 
Following the work of Ding and Tian it remained an open problem how to extend their techniques
and provide a complex-geometric proof of the Moser-Trudinger-Onofri inequality.
The key hurdle in proving \Ref{MTOThm} is to extend the argument to the set
$\CinfM\setminus\calH_\o$ which a priori has no clear geometric significance as it 
represents indefinite forms rather than \K forms.

Alternatively, what is missing is a geometric interpretation of the Berger-Moser-Ding functional
$F_1$.
The following result is the key ingredient in our proof of \MTClassicalineq\
\ref{R2, Lemma 2.4}. 
\FProp{
\TaggProp{EnFProp}
Let $\O=c_1$. The following relation holds
$$
(\Ric^{-1})^\star E_n=F_1,\q \hbox{on \ } \Hc\times\Dc.
$$
}

This provides a geometric interpretation for $F_1$. Indeed, the functional $E_n$ is 
the potential for the Laplacian of the determinant of the Ricci tensor, considered as a 1-form on $\HO$,
i.e., $dE_n(\o,\ovp)=\D_{\ovp}\Big(\frac{(\Ric\ovp)^n}{\ovpn}\Big) \ovpn$.

Hence, this result explains the geometric meaning the set $\CinfM\setminus\Ho$ plays in
the Moser-Trudinger-Onofri inequality. Namely, a function will satisfy this inequality
if and only if it represents the Ricci form of a \K metric whose Ricci energy $E_n$ is nonnegative
with respect to a \KE metric. It now becomes important to understand the sets
$\calA_n=\{\ovp\in\Hc\,:\, E_n(\o,\ovp)\ge0\}$, defined in \Ref{SectionEnergy}.
Naturally, we introduce the following definition.

\FDef{
Let the Moser-Trudinger-Onofri neighborhood of $\Ho$ be the subset
$$
MTO_n=\{\vp\in\CinfM: \!\vp \hbox{$\!$ satisfies \GeneralizedMTOeq\ on the Fano manifold\ } 
(M,\JJJ), \dim_\CCfoot M=n\}.
\eqn\MTOnEq
$$ 
}

We are now in a position to state our generalization of \Ref{DTBoundednessThm}
that is optimal in higher dimensions as well as some practical bounds. 
The result says that the Moser-Trudinger-Onofri inequality holds in higher dimensions on a canonically defined
set $MTO_n$ that is strictly larger than the space of \K potentials $\Ho$ and is geometrically 
related to Ricci curvature.

\FThm{
\TaggThm{CharacterizationThm}
Let $(M,\JJJ,\o)$ be a \KE Fano manifold.\newline
(i) The generalized Moser-Trudinger-Onofri inequality \GeneralizedMTOeq\ holds precisely 
on the set $MTO_n=\Ric(\calA_n)$. Furthermore, $\Hc\subsetnoteq MTO_n\sseq \Dc$.\newline
(ii) Define the sets $\calB_k:=\{\ovp\in\Hc\,:\, I_k(\ovp,\Ricovp)\ge0\}$. Then one
has $\Hcplus\subsetnoteq\calB_k\sseq\calA_k$.\newline
(iii) One has $\calA_1=\calB_1=\Hc$, $\calB_2\supseteq\{\ovp\in\Hc\,:\, \Ricovp+2\ovp\ge0\}$,
$\calB_3\supseteq\{\ovp\in\Hc\,:\, \Ricovp+\ovp\ge0\}$, and for each $k$ one
may readily obtain an explicit bound on the set $\calB_k$, hence on $\calA_k$, in terms
of a lower bound on the Ricci curvature using \Ikeq. In particular there exist $c_n>0$ depending only on $n$ such that,
$$
MTO_n\supseteq\{\vp\in\CinfM\,:\, \ovp\ge -c_n\ricm1\ovp\},
$$
and, e.g., $c_1=\infty, c_2\ge2, c_3\ge1$.
}

As a corollary we are now able to provide a complex-geometric proof of the classical 
Moser-Trudinger-Onofri inequality.

\bigskip
\noindent
{\it Proof of \Ref{MTOThm}.} Observe that $MTO_1=\Ric(\calA_1)=\Ric(\Hc)=\Dc$. The last
equality requires solving the equation%
\note{\cmreight This is the classical $\meight n=1$ version of the Calabi-Yau Theorem whose proof goes back at least
to\hfill\break
\vglue-0.5cm
\hglue-\parindent\cmreight
Wallach and Warner \reffoot{WW}.}
$\Ric\o_\vp=\o_\psi$ for $\vp$, equivalently Poisson's equation $\D_\o\vp=e^{\fo-\psi}-1$.\done

\subsubsection{An analytic characterization of \KE manifolds and an analytic criterion for almost-\KE manifolds}
\TaggSubS{AnalyticCharacterizationSubS}%
In the first part of this subsection we explain how the inverse Ricci operator can be used to solve a
problem concerning energy functionals on the space of \K forms. 
We hope this sheds new light on the solution of this problem that we gave previously \ref{R2}.

{\TaggPage{AnalyticCharacterizationPage}}%

Chen and Tian's generalization of Mabuchi's \K energy, $E_0$, and of Bando and Mabuchi's 
Ricci energy, $E_n$, to a family
of functionals $\{E_k\}_{k=0}^n$ (see \Ref{SectionEnergy} for definitions) naturally 
raised the question of whether Tian's analytic characterization of \KE manifolds
in terms of $E_0$ generalizes to these functionals. 
In addition it raised the question whether Bando and Mabuchi's 
criterion for almost-\KE manifolds in terms of $E_0$ generalizes to these functionals.
These questions were also independently raised by Chen \ref{Che2, p. 37; CLW, \S1.3}.
We now recall both of these fundamental results and explain how to generalize them.
This provides an answer to these questions. It shows that the answer is both ``yes" and 
``no": these criteria extend to the other functionals $\{E_k\}$, however they fail to extend 
in an identical manner. The subtlety comes from the appearance of the inverse Ricci
operator as we will see below.

\FThm{
\TaggThm{BMDTThm}
Let $(M,\JJJ)$ be a Fano manifold.\newline
(i) \ref{B,BM,DT} If either $F_1$ or $E_0$ is bounded from below on $\Hc$ then for 
each $\eps>0$ there exists a \K metric $\o_\eps\in\Hc$ satisfying
$\Ric\o_\eps>(1-\eps)\o_\eps$.
\newline
(ii) \ref{T5,T6,TZ1} Assume that $\AutMJ$ is finite.%
\note{\cmreight In the general case a slightly more involved statement holds (see \reffoot{R2} for details).} 
Then the properness of 
$F_1$ (or $E_0$) on $\Hc$ is equivalent to the existence of a \KE metric.
}
Our strategy in extending these results to the functionals $\{E_k\}_{k=0}^n$ was: 
(a) first prove a new formula that expresses $E_k$ in terms of the sum of $E_0$ 
and another new exact energy functional 
$(\o,\ovp)\mapsto I_k(\ovp,\Ric\ovp)-I_k(\o,\Ric\o)$ (\Ref{BMRProp}) and use it to show
$$
\eqa{
F_1 \hbox{\ bounded from below on\ } \Hc 
& \Ra
E_0 \hbox{\ bounded from below on\ } \Hc
\cr & \Ra
E_1 \hbox{\ bounded from below on\ } \Hc
\cr & \Ra
E_2 \hbox{\ bounded from below on\ } \Hcplus
\cr & \;\;\vdots
\cr & \Ra 
E_n \hbox{\ bounded from below on\ } \Hcplus.
}
$$ 
(b) Next use \Ref{EnFProp} to conclude:
$$
E_n \hbox{\ bounded from below on\ } \Hcplus
\Ra
F_1 \hbox{\ bounded from below on\ } \Hc.
$$
(c) Finally, some additional arguments were needed in order to prove that the 
properness of $E_n$ on $\Hcplus$ implies the existence of a \KE metric. 

We can now state the extension of the 
theorems of Bando-Mabuchi and Tian to the energy functionals $\{E_k\}$. The case $k=1$ 
was proven before by Chen-Li-Wang and Song-Weinkove in a different manner  \ref{CLW,SW}.

\FThm{
\TaggThm{BoundednessPropernessThm}
Let $(M,\JJJ)$ be a Fano manifold.\newline
(i) If either $F_1$ or $E_k$ (for some $k\in\{0,\ldots,n\}$) is bounded from below on 
$\Hcplus$ then for each $\eps>0$ there exists a \K metric $\o_\eps\in\Hc$ satisfying
$\Ric\o_\eps>(1-\eps)\o_\eps$.
\newline
(ii) Assume that $\AutMJ$ is finite.\note{\cmreight See footnote to \Ref{BMDTThm}.}
Then the properness of 
$F_1$ or of $E_k$ (for some $k\in\{0,\ldots,n\}$) on $\Hcplus$ is 
equivalent to the existence of a \KE metric.
}

It is important to note that the appearance of the inverse Ricci operator in step
(b) was a crucial ingredient. The discrepancy between the behavior of $F_1,E_0,E_1$ 
and that of the
functionals $E_2,\ldots,E_n$ can be explained using the time one Ricci iteration: the first three 
are unconditionally monotone along the iteration, while for the latter $n-1$ this is 
true if and only if one assumes that the initial point lies in $\calB_k$,
 and
$\Hcplus\subsetnoteq\calB_k\sseq\calA_k\subsetnoteq\Hc$ 
(\Ref{EnergyMonotonicityProp}). 
Furthermore, along the first step of the iteration the functionals $E_k$ may increase by an arbitrary amount!
To be precise, for any Fano manifold (\KE or not) we have the following result \ref{R2}:\note{%
\cmreight We believe that the same result should hold for $\meight E_k$ for each $\meight 2\le k\le n$.}
The Ricci energy $E_n$ is bounded from below on $\Hc$ if and only if $n=1$.
We conclude that the assumption in \Ref{BoundednessPropernessThm} (ii) %and in \EkLimitEq\ 
is essential and cannot be weakened
from $\Hcplus$ to $\Hc$. This explains our remark earlier on the subtlety present when $k\ge2$.

Previously, several authors (for references see \ref{R2}) have proven the four implications
on the equivalence of the boundedness from below of $F_1$, $E_0$ and $E_1$. 
Often they appealed to results on
the Ricci flow. This suggests that in this context the Ricci iteration rather
than the flow is perhaps more suited.

\FRemm{In light of the discussion above it would be interesting to know whether
there exist intial conditions in $\Hc$ for which $E_k, \;k\ge2$, increases by an arbitrary amount
along the Ricci flow restricted to the time interval $[0,1]$.
}

\subsubsection{A new Moser-Trudinger-Onofri inequality on the Riemann sphere 
and a family of energy functionals}
\TaggSubS{ImprovedMTOEnergyFunctionalsSubS}%
{\TaggPage{ImprovedMTOEnergyFunctionalsPage}%
In the first part of this subsection we prove 
results that improve on the restricted generalized Moser-Trudinger-Onofri inequality 
(\Ref{DTBoundednessThm})
in a different direction than that explored in \Ref{MTOSubS} (\Ref{CharacterizationThm}).
Namely, we show that the inequality holds even when one adds certain negative terms to the exponent
on the right hand side. This is done by expressing the excess in the inequality in geometric terms, 
namely in terms of the inverse Ricci operator.
This is different from Tian's approach to a strenghtened inequality on \KE manifolds \ref{T6, Theorem 6.21} 
and in particular involves sharp constants and a precise characterization of the case of equality. 
In the future we hope to address the relation between these two approaches.
In the second part we introduce a family of energy functionals and explain their relation
to the improved inequality.

We now state the main result of this subsection. 

\FThm{
\TaggThm{MTOFirstIneqThm}
Let $(M,\JJJ,\o)$ be a Fano \KE manifold.
Then for each $\vp\in MTO_n$ holds
$$
\eqa{
\V\intm e^{-(\vp-\V\intm\vp)}\on
\le
e^{
J(\o,\ovp)-{\h{$\mf\sum$}}_{j=1}^\infty J(\Ricnotsosmall^{-j}\ovp,\Ricnotsosmall^{-j+1}\ovp)
}.
}\eqn\ImprovedMTOIneq$$
Each of the terms in the sum is nonnegative precisely when 
$\ovp\in\Ric(\calB_n)\supsetnoteq\Hc$.%
\note{\cmreight More precisely, each of the terms with $\meight j\ge2$ is nonnegative for all 
$\meight \vp\in MTO_n$, while the term\break
\vglue-0.5cm
\hglue-\parindent\cmreight with $\meight j=1$ is nonnegative precisely
when $\meight \ovp\in\h{\cmreight Ric}\,(\calfootB_n\!)$. Therefore, after possibly omitting the first\break
\vglue-0.5cm
\hglue-\parindent\cmreight term in
the sum, \ImprovedMTOIneq\ is an improvement over \GeneralizedMTOeq\ for all
$\meight \vp\in MTO_n$ and not just for the subset\break
\vglue-0.5cm
\hglue-\parindent\cmreight $\meight\h{\cmreight Ric}\,(\calfootB_n\!)\sseq MTO_n$ (both
strictly contain $\meight\Hcfoot$).}
Equality holds if and only if $\ovp$ is the pull-back of $\o$ by a holomorphic transformation.
}

\noindent

Recall that \Ref{CharacterizationThm} strengthened the 
restricted generalized Moser-Trudinger-Onofri inequality
(\Ref{DTBoundednessThm}) by optimally enlarging the set of functions on which it holds
to a set strictly containing $\Ho$. \Ref{MTOFirstIneqThm} further strengthens
\Ref{CharacterizationThm}: it shows that the sets $MTO_n$ (see \MTOnEq) 
are characterized by an inequality stronger than \GeneralizedMTOeq.
A version of this result holds also under the assumption that the K-energy is bounded from below.
For simplicity we only state the result in the \KE setting.

\medskip
\proof
By definition $F_1(\o,\ovp)\ge 0$ for each $\vp\in MTO_n$. Observe that by \Ref{CharacterizationThm} (i)
it follows that $\Ric^{-1}$ preserves
$MTO_n$. Therefore for each $l\in\NN$,
$$
F_1(\o,\Ric^{-l} \ovp)\ge 0, \q\all\,\vp\in MTO_n.
$$
By exactness of $F$ we obtain
$$
F_1(\o,\ovp)+F_1(\ovp,\Ric^{-1} \ovp)+\ldots + F_1(\ricm {l+1} \ovp,\ricm l \ovp)\ge 0,
$$
that is,
$$
F_1(\o,\ovp)\ge \sum_{j=1}^{l}  F_1(\ricm {j} \ovp,\ricm{j+1} \ovp).\eqn\FIteratedIneq
$$
Now, using \Ieq-\Jeq\ and \FIterationSecondEq\ one has
$$
F_1(\o,\ovp)
=
J(\o,\ovp)-\V\intm\vp\on-\log\V\intm e^{f_\o-\vp}\on. \eqn\FFunctionalSecondEq
$$
It follows that for any $\a\in\Hc$ holds
$$
F_1(\alpha,\Ric\alpha)=J(\alpha,\Ric\alpha)-\V\intm f_\a\alpha^n.\eqn\FIteratedSecondEq
$$
Combining \FIteratedIneq-\FIteratedSecondEq, and letting $l$ tend to infinity, 
yields
$$
\eqa{
\V\intm e^{-(\vp-\V\intm\vp)}\on
& \le
e^{
J(\o,\ovp)-{\h{$\mf\sum$}}_{j=1}^\infty J(\Ricnotsosmall^{-j}\ovp,\Ricnotsosmall^{-j+1}\ovp)
}
\cr & \qq
\cdot e^{
{\h{$\mf\sum$}}_{j=1}^\infty\V\intm f^{(j)}_{\ovp}(\Ricnotsosmall^{-j} \ovp)^n,
}
}\eqn\ComplicatedMTOJIneq
$$
where $f^{(j)}$ is the push-forward of the vector field $f$ under $\ricm j$.
Since by \RicciPotentialEq\ the second term in \FIteratedSecondEq\ is nonnegative the desired inequality
now follows from \ComplicatedMTOJIneq. 

The last statement follows from the fact that $I_n=J$, \Ref{EnFProp}, and the definition of $\calB_n$ 
\calBkEq.
\done

In the case of the Riemann sphere $S^2$, \Ref{CharacterizationThm} (iii) implies 
$\Ric(\calB_1)=\Ric(\calA_1)=MTO_1=\Cinf(S^2)$. 
Therefore we have the following improvement of the 
classical Moser-Trudinger-Onofri inequality (\Ref{MTOThm}). For notation we refer to \Ref{MTOSubS}.

\FCor{
\TaggCor{MTOSphereCor}
Denote by $(S^2,\JJJ,\o=\o_{\hbox{\small FS},2/V})$ a round sphere of volume $V$. 
For any function $\vp$ on $S^2$ in $W^{1,2}(S^2)$ one has
$$
\V\int_{S^2} e^{-\vp+\V\int_{S^2}\vp\o}\o
\le 
e^{
\V\int_{S^2}\ha\i\del\vp\w\dbar\vp
-{\h{$\mf\sum$}}_{j=1}^\infty J(\Ricnotsosmall^{(-j)}\ovp,\Ricnotsosmall^{(-j+1)}\ovp)
}.\eqn\ImprovedClassicalMTOIneq
$$
Each of the terms in the sum is nonnegative, and equals zero if and only if
$\ovp$ is obtained from $\o$ by a M\"obius transformation. This also characterizes
when equality holds in \ImprovedClassicalMTOIneq.
}

Note that the smoothing property of the iteration (see page \Ref{SmoothningPropertyPage})
implies that the extra terms in the sum are meaningful under the assumption $\vp\in W^{1,2}(S^2)$.

Motivated by \Ref{EnFProp} we define the following family of energy functionals.
For each $k\in\{0,\ldots,n\}$ and $l\in\NN\cup\{0\}$ let $E_{k,l}$ denote the pull-back
by $\ricm l$ of the Chen-Tian functional $E_k$ (see \ChenTianFunctionalsEq).
That is
$$
E_{k,l}(\o,\ovp)=E_k(\Ric^{-l} \o,\Ric^{-l}\ovp).\eqn\ElkFunctionalsEq
$$
For example, $E_{n,1}=F_1$, and
$$
E_{k,1}(\o,\ovp)=
F_1(\o,\ovp)-(J-I_{k})(\ricm1\ovp,\ovp)+(J-I_{k})(\ricm1\o,\o),\q k=0,\ldots,n.\eqn
$$

In light of this, \Ref{MTOFirstIneqThm} is seen to be a corollary of the following inequality:
$$
E_{n,l+1}(\o,\;\cdot\;)=(\ricm l)^\star F_1\,\big|_{\{\o\}\times MTO_n}\ge0, \q \all\, l\in\NN.\eqn
$$

Finally, we remark that in light of \Ref{GeneralizedInverseRicciDef} and \Ref{EnFProp} 
one may also extend the definition of Ding's functional to an arbitrary \K manifold and class. This
might have some future applications.

\subsubsection{Construction of Nadel-type obstruction sheaves}
\TaggSubS{NadelSheavesSubS}%
{\TaggPage{NadelSheavesPage}}%
Up until this point we have scarcely concerned ourselves with the behavior of 
the various dynamical systems constructed in the absence of a fixed point.
In this subsection we show that in this situation, and in the Fano setting, 
the Ricci iteration will produce Nadel-type obstruction sheaves, similarly to the
continuity method and the Ricci flow. 
The basic references for this subsection are Demailly-Koll\'ar \ref{DK} and Nadel \ref{N1}.

Let $PSH(M,\JJJ,\o)\sseq\Loneloc(M)$ denote the set of $\o$-plurisubharmonic functions.
For $\vp\in PSH(M,\JJJ,\o)$ define the multiplier ideal sheaf associated to $\vp$ as the
sheaf $\calI(\vp)$ defined for each open set $U\sseq M$ by local sections
$$
\calI(\vp)(U)=\{h\in \calO_M(U): |h|^2 e^{-\vp}\in\Loneloc(M) \}.\eqn\IdealSheafSectionsEq
$$
Such sheaves are coherent. % \ref{D, p. 73; N1}. 
Such a sheaf is called proper if it is neither zero nor the structure sheaf $\calO_M$.

Nadel showed that in the absence of a \KE metric the continuity method \AubinPathEq\ will
produce a certain family of multiplier ideal sheaves.
Phong, \v Se\v sum and Sturm showed that certain multiplier ideal sheaves can be obtained
also from the Ricci flow
$$
\ovptn=\on e^{\fo-\vpt+\dot\vpt},\q \vp(0)=\h{const}.\eqn\RicciFlowPositiveMAEq
$$
Denote by $\lfloor x\rfloor$ the largest integer not larger than $x$.

\FThm{
\TaggThm{PSSThm}
\ref{PSS}
Let $(M,\JJJ)$ be a Fano manifold not admitting a \KE metric. 
Let $\gamma\in(1,\infty)$ and let $\o\in\Hc$.
Then there exists an initial condition $\vp(0)$ and a subsequence 
$\{\vp_{t_j}\}_{j\ge0}$ of solutions of 
\RicciFlowPositiveMAEq\
such that $\lim_{j\ra\infty}\vp_{t_j}=\vp_\infty\in PSH(M,\JJJ,\o)$ and
$\calI(\gamma\vp_\infty)$ is a proper multiplier ideal sheaf satisfying 
$$
H^r(M,\calI(\gamma\vp_\infty)\otimes K_M^{-\lfloor\gamma\rfloor})=0, \q \all\, r\ge 1.
\eqn\NadelVanishingFanoEq
$$
}

Their proof relies on some of Perelman's estimates for the Ricci flow as well as the following theorem 
of \Kolodziej.

\FThm{\TaggThm{KolodziejThm}
\ref{Ko}
Let $F\in L^p(M,\o), p>1$ be a positive continuous function with $\V\intm F\on=1$. 
There exists a bounded solution $\vp$ to the equation $\ovpn= F\on$ on $M$ 
which satisfies $\osc\vp\le C$ with $C$ depending only on 
$||F||_{L^p(M,\o)},p$ and $(M,\o)$.
}

Let $\tau=1/\mu=1$. The following simple result is a discrete analogue of \Ref{PSSThm}. Its
very simple proof compared to that of the analogous result for the Ricci flow is our 
main motivation for including it here. Moreover, the sheaves produced in this way are essentially
computable (see the next subsection).

\FThm{
Let $(M,\JJJ)$ be a Fano manifold not admitting a \KE metric. 
Let $\gamma\in(1,\infty)$ and let $\o\in\Hc$.
Then there exists a subsequence 
$\{\psi_{j_k}\}_{k\ge1}$ of solutions of 
\IterationKEEq\
such that $\lim_{k\ra\infty}\psi_{j_k}=\psi_\infty\in PSH(M,\JJJ,\o)$ and
$\calI(\gamma\psi_\infty)$ is a proper multiplier ideal sheaf satisfying 
\NadelVanishingFanoEq.
}

\proof 
Indeed, since the iteration takes the form
$$
\o_{\psi_{l+1}}^n=\on e^{f_\o-\psi_{l}},\q l\in\NN,
$$
\Ref{KolodziejThm} can be directly 
applied (observe that from \RicciPotentialEq\ an estimate on $\osc\psi_l$ implies one on
$||\psi_l||_{\Linf(M)}$) to construct sheaves with $\gamma>1$, making use of 
\Ref{EnergyMonotonicityProp} 
(for more details see \ref{R4, \S2 (iv)}).\done

\FRemm{One may also construct multiplier ideal sheaves for the Ricci iteration with other time steps and for
exponents in the range $(n/(n+1),1)$ much the same as the continuity method sheaves constructed by Nadel
as well as the analogous ones constructed in \ref{R4} for the Ricci flow
(we hope to discuss this in more detail in the sequel; note that the latter construction strenghthened \Ref{PSSThm}).
}

We remark that in the context of this section, it is also interesting to study the limiting behavior of the 
inverse Ricci operator (\Ref{AnotherFlowSection}) under iteration.

\subsubsection{Relation to balanced metrics}
\TaggSubS{BalancedMetricsSubS}%
{\TaggPage{BalancedMetricsPage}}%
In this paragraph we describe an immediate corollary of the work of Donaldson. 
It was pointed out to me by J. Keller.
It gives with no further work an algorithm for computing \KE metrics using balanced metrics: 
Given a polarized Hodge manifold $(X,L)$ and a volume form $\nu$ 
Donaldson \ref{Do3} constructs a sequence of pull-backs of Fubini-Study metrics in $\calH_{c_1(L)}$ 
induced from Kodaira embeddings that converge to a solution of the Calabi-Yau equation
$\o_\vp^n=\nu$. Since in the Fano case our time one Ricci iteration consists precisely of solving a Calabi-Yau
equation at each iteration we see that repeated application of 
Donaldson's constructions approximates the Ricci iteration and in this sense provides a quantization 
of the Ricci flow.

Another consequence is the possibility to numerically construct Nadel-type sheaves
on Fano manifolds admitting no \KE metrics, by \Ref{NadelSheavesSubS}.

Note that more generally one may approximate in the same manner the orbits
of the iteration given by the inverse Ricci operator (\Ref{GeneralizedInverseRicciDef}) on
an arbitrary \K manifold with $\O\in H^2(M,\ZZ)$.

Finally, it would be interesting to find more relations between discretizations of other
geometric flows and iteration schemes involving Bergman metrics.

\subsubsection{A question of Nadel}
\TaggSubsection{NadelQuestionSubsection}
\TaggSubS{NadelQuestionSubS}
{\TaggPage{NadelQuestionPage}}
\unskip\unskip As explained in \Ref{CanonicalMetricsSection} one of the original motivations 
for our work was a question raised by Nadel \ref{N2}: Given 
$\o\in\Hc$ define a sequence of metrics $\o,\Ric\o,\Ric(\Ric\o),\ldots,$ as long as positivity 
is preserved; what are the periodic orbits of this dynamical system? The cases $k=2,3$ 
in the following theorem are due to Nadel. 

\FThm
{\TaggThm{NadelIterationThm}
Let $(M,\JJJ,\o)$ be a Fano manifold and assume that $\Ric^{l}\,\o=\o$
for some $l\in\ZZ$. Then $\o$ is \KEno.
}

\proof
The theorem follows from \Ref{EnergyMonotonicityProp}. Indeed,
note that the nonexistence of periodic fixed points of negative order implies that
of positive order, and vice versa. 
Therefore assume that for some $\o\in\Hc$ and some
$l\in\NN$ one has $\ricm l\,\o=\o$. By the cocycle condition we thus have
$$
0=E_0(\o,\Ric^{-l}\o)=\sum_{i=0}^{l-1} E_0(\Ric^{-i}\o,\Ric^{-i-1}\o).\eqn\EzeroCocyleEq
$$
By \Ref{EnergyMonotonicityProp} one has
$$
E_0(\Ric^{-i}\o,\Ric^{-i-1}\o)<0,
$$ 
unless $\ricm{i}\o=\ricm {i-1}\o$. 
Therefore each of the terms in \EzeroCocyleEq\ vanishes and $\o$ is \KEno.\done
Moreover, from the proof we have the following stronger conclusion:
\FCor{
\TaggCor{NadelIterationCor}
Let $(M,\JJJ,\o)$ be a Fano manifold with trivial Futaki character.
Assume that $\Ric^{l}\,\o=h^\star\o$
for some $l\in\ZZ$ and some $h\in\AutMJ$. Then $h=\id$ and $\o$ is \KEno.
}

\Ref{TwistedEnergyMonotonicityLemm} implies the following natural generalization of
\Ref{NadelIterationThm} to the setting of solitons.

\FCor{
Let $(M,\JJJ,\o)$ be a Fano manifold and assume that $X$ 
belongs to a reductive Lie subalgebra of $\autMJ$ and that the 
one-parameter subgroup $T_{\h{\notsosmall J} X}$ generated by $\JJJ X$ is 
a compact torus in $\AutMJ$. Let $\o\in\Hc(T_{\h{\notsosmall J} X})$. Assume that 
$\Ric_{\psi^X}^{l}\,\o=\o$
for some $l\in\ZZ$. Then $\o$ is a \KR soliton.
}

In addition, under the assumption that the Tian-Zhu character \ref{TZ3} is trivial one has a statement
analogous to \Ref{NadelIterationCor}. Also, as noted in \Ref{RicciIterationTwistedSection}, and using
a generalized character introduced by Futaki \ref{F2}, this result
extends to the setting of multiplier Hermitian structures.

To conclude this subsection we remark that what now becomes apparent is that 
Nadel's iteration scheme is precisely the Euler method for the conjugate Ricci flow
and is thus dual to our iteration that corresponds to the backwards Euler method for
the Ricci flow. %This completely elementary observation lies at the basis of this article.

\FRemm{
In light of \Ref{NadelIterationThm} perhaps it would be interesting to re-examine Nadel's
generalized maximum principle which was used to provide a completely different proof for the
cases $k=2,3$.
}

\subsubsection{The Ricci index and a canonical nested structure on the space of \K metrics}
{\TaggPage{RicciIndexPage}}%
In this subsection we describe a new canonical structure inherent in the space of \K 
forms determined by the complex structure and the \K class alone.

Consider first the case of a Fano manifold.
As we saw earlier the iteration of the inverse Ricci operator on $\Hc$ has the advantage
of possessing infinite orbits starting at any initial points. The Ricci operator on the other
hand lacks this property, according to the Calabi-Yau theorem. This motivates the following
definition.

\FDef{
\TaggDef{NestedStructureDef}
Let $(M,\JJJ)$ be a Fano manifold. For each $l\in\NN\cup\{0\}$ 
denote by $\hc l$ denote the domain of definition of $\Ric^{l}$. 
}
One has 
$$
\Dc=\hc 0\supset \Hc=\hc 1\supset\hc 2=\Hcplus\supset\cdots\supset\hc l\supset\cdots.\eqn
$$
In other words, we may define on $\Hc$ an integer-valued function 
$$
\o\mapsto r(\o),\eqn
$$
where $r(\o)$ is the unique positive integer satisfying 
$\o\in\hc{r(\o)}\setminus \hc{r(\o)+1}$. When no such number exists we set
$r(\o)=\infty$. We call the function $r:\Hc\ra\NN$ the Ricci index.
The number $r(\o)$ is a Riemannian invariant of the manifold $(M,\JJJ,\o)$.
It may also be defined for general Riemannian manifolds however it seems hard
to study in such generality.

One may extend such a construction to a general \K manifold in at least two ways, using either
the Ricci iteration or the inverse Ricci operator. Choosing the latter we obtain
the following extension of \Ref{NestedStructureDef}.

\FDef{
\TaggDef{GeneralizedNestedStructureDef}
Let $(M,\JJJ)$ be a \K manifold and let $\O$ denote a \K class. 
For each $l\in\NN$ denote by $\ho l$ the image of $\HO$
under $\Ric_\O^{-l+1}$. 
}

Several natural questions arise that we hope to touch upon in the future. 
What is $\ho\infty:=\bigcap_{l=1}^\infty\ho l$? How
to asymptotically relate the Ricci index to the time parameter of the Ricci flow?
Also, how to relate the Ricci index, on the one hand, to the metric structure on the 
space $\HO$ \ref{M2,Se,Do2} defined by
$\langle\mu,\nu\rangle_\o=\V\intm\mu\nu\on,\;\all\mu,\nu\in T_\o\HO\isom\CinfM/\RR$ 
and, on the other hand, to sublevel sets of Calabi's energy and Mabuchi's K-energy? Finally, what is the
relation between the Ricci index and positivity?

\def\smallblackbox{\vrule height.6ex width .5ex depth -.1ex}
\def\boxseparation{\hfil\smallblackbox$\q$\smallblackbox$\q$\smallblackbox\hfil}
\bigskip
\boxseparation
\bigskip

I would like to express my deep gratitude to my teacher, Gang Tian, for 
his advice, warm encouragement and for pointing out to me the relevance of \ref{TY}. 
I am indebted to him as well as to J. Song for suggesting
that the inverse Ricci operator I defined could be related to the Ricci flow. This had
a decisive impact on the present work.
I would also like to thank Curtis McMullen, whose class several years ago has
been a source of great inspiration.
The first version of this article \ref{R1} was written while I was a Visiting Scholar
at Peking University during Summer 2005.  The present version was 
completed in Summer 2007 during a visit to the Technion. I thank
both institutions for their hospitality and partial financial support.
Much of the present work was presented in the past year at several places, including 
at Imperial College in December 2006.
This material is based upon work supported under a National Science 
Foundation Graduate Research Fellowship.

\frenchspacing

\bigskip\bigskip
\noindent{\bf Bibliography}
\TaggPage{BibliographyPage}

\bigskip
\def\ref#1{\Taggf{#1}\item{ {\bf[}{\sans #1}{\bf]} } }

\def\sm{\vglue1pt}

\ref{A1} Thierry Aubin, \'{E}quations du type {M}onge-{A}mp\`ere sur les vari\'et\'es
k\"ahl\'eriennes compactes, {\sl Bulletin des Sciences Math\'ematiques} 
{\bf 102} (1978), 63--95.
\sm
\ref{A2} \opcit, R\'eduction du cas positif de l'\'equation de {M}onge-{A}mp\`ere 
sur les vari\'et\'es k\"ahl\'eriennes compactes \`a la d\'emonstration d'une in\'egalit\'e, 
{\sl Journal of Functional Analysis} {\bf 57} (1984), 143--153.
\sm
\ref{A3} \opcit, Some nonlinear problems in Riemannian Geometry, Springer, 1998.
\sm
\ref{BE} Dominique Bakry, Michel {\'E}mery, Diffusions hypercontractives, in
{\it S\'eminaire de probabilit\'es XIX, 1983/84} (J. Az\'ema et al., Eds.),
Lecture Notes in Mathematics {\bf 1123}, Springer, 1985, 177--206.
\sm
\ref{B} Shigetoshi Bando, The K-Energy Map, Almost \KE Metrics and an
Inequality of the Miyaoka-Yau Type, {\sl T\^ohoku Mathematical Journal}
{\bf 39} (1987), 231--235.
\sm
\ref{BM} Shigetoshi Bando, Toshiki Mabuchi, Uniqueness of \KE metrics
modulo connected group actions, in {\it Algebraic Geometry,
Sendai, 1985} (T. Oda, Ed.), Advanced Studies in Pure Mathematics {\bf 10},
Kinokuniya, 1987, 11--40.
\sm
\ref{Be} William Beckner, Sharp Sobolev inequalities on the sphere and the 
Moser-Trudinger inequality, {\sl Annals of Mathematics} {\bf 138} (1993), 213--242.
\sm
\ref{Ber} Melvyn S. Berger, Riemannian structures of prescribed Gaussian curvature for
compact 2-manifolds, {\sl Journal of Differential Geometry} {\bf 5}
(1971), 325--332.
\sm
\ref{Bes} Alfred L. Besse, Einstein manifolds, Springer, 1987.
\sm
\ref{C1} Eugenio Calabi, The variation of \K metrics. I. The structure of the space; II. A minimum problem, 
{\sl Bulletin of the American Mathematical Society} {\bf 60} (1954), 167--168.
\sm
\ref{C2} \opcit, On K\"ahler manifolds with vanishing canonical class,
in {\it Algebraic geometry and topology. A symposium in honor of S. Lefschetz} (R. H. Fox, Ed.), 
Princeton University Press, 1957, 78--89.
\sm
\ref{Ca} Huai-Dong Cao, Deformations of \K metrics to \KE metrics on compact
\K manifolds, {\sl Inventiones Mathematicae} {\bf 81} (1985), 359--372.
\sm
\ref{CL1} \opcit, Competing symmetries of some functionals arising in mathematical physics, 
in {\it Stochastic processes, physics and geometry} (S. Albeverio et al., Eds.), World
Scientific, 1990, 277--288.
\sm
\ref{CL2} Eric A. Carlen, Michael Loss, Competing symmetries, the logarithmic HLS inequality and Onofri's
inequality on $S^n$, {\sl Geometric and Functional Analysis} {\bf 2} (1992), 90--104.
\sm
\ref{Ch} Sun-Yung A. Chang, Non-linear elliptic equations in conformal geometry,
European Mathematical Society, 2004.
\sm
\ref{Che1} Xiu-Xiong Chen, On the lower bound of the Mabuchi energy and its application,
{\sl International Mathematics Research Notices} (2000), 607--623.
\sm
\ref{Che2} \opcit, On the lower bound of energy functional $E_1$ (I)---a stability
theorem on the K\"ahler-Ricci flow, {\sl The Journal of Geometric Analysis} {\bf 16}
(2006), 23--38.
\sm
\ref{CLW} Xiu-Xiong Chen, Hao-Zhao Li, Bing Wang,
On the K\"ahler-Ricci flow with small initial $E_1$ energy (I),
 preprint, arxiv: math.DG/0609694 v2. To appear in {\sl Geometric and Functional Analysis}.
\sm
\ref{CT} Xiu-Xiong Chen, Gang Tian, Ricci flow on \KE surfaces,
{\sl Inventiones Mathematicae} {\bf 147} (2002), 487--544.
\sm
\ref{Cho} Bennett Chow et al., The Ricci flow: Techniques and applications. Part I: Geometric
aspects, American Mathematical Society, 2007.
\sm
\ref{DK} Jean-Pierre Demailly, J\'anos Koll\'ar, Semi-continuity of complex singularity exponents
and K\"ahler-Einstein metrics on Fano orbifolds, {\sl Annales Scientifiques de l'\'Ecole
Normale sup\'erieure} {\bf 34} (2001), 525--556.
\sm
\ref{D} Wei-Yue Ding, Remarks on the existence problem of positive
{K}\"ahler-{E}instein metrics, {\sl Mathematische Annalen} {\bf 282} (1988), 463--471.
\sm
\ref{DT} Wei-Yue Ding, Gang Tian, The generalized Moser-Trudinger inequality, 
in {\it Nonlinear Analysis and Microlocal Analysis: Proceedings of the International
Conference at Nankai Institute of Mathematics} (K.-C. Chang et al., Eds.), 
World Scientific, 1992, 57--70. ISBN 9810209134.
\sm
\ref{Do1} Simon K. Donaldson, Anti self-dual Yang-Mills connections over complex
algebraic surfaces and stable vector bundles, {\sl Proceedings
of the London Mathematical Society} {\bf 50} (1985), 1--26.
\sm
\ref{Do2} \opcit, Symmetric spaces, K\"ahler geometry and Hamiltonian dynamics, in
{\it Northern California Symplectic Geometry Seminar} (Ya. Eliashberg et al., Eds.),
American Mathematical Society Translations: Series 2 {\bf 196},
American Mathematical Society, 1999, 13--33.
\sm
\ref{Do3} \opcit, Some numerical results in complex differential geometry, preprint, 
April \thhnotsosmall{27}, 2006.
\sm
\ref{F1} Akito Futaki, K\"ahler-{E}instein metrics and integral invariants,
Lecture Notes in Mathematics {\bf 1314}, Springer, 1988.
\sm
\ref{F2} \opcit, Some invariant and equivariant cohomology classes of the space of {K}\"ahler metrics,
{\sl Proceedings of the Japan Academy. Series A} {\bf 78} (2002), 27--29.
\sm
\ref{Gh} Alessandro Ghigi, On the Moser-Onofri and Pr\'ekopa-Leindler inequalities,
{\sl Collectanea Mathematica} {\bf 56} (2005), 143--156.
\sm
\ref{G1} Daniel Z.-D. Guan, Quasi-Einstein metrics, {\sl International Journal of Mathematics}
{\bf 6} (1995), 371--379.
\sm
\ref{G2} \opcit, Extremal-solitons and $\Cinf$ convergence of the modified Calabi flow on certain
$CP^1$ bundles, preprint, December 2\nd\unskip, 2006.
\sm
\ref{H} Richard S. Hamilton, Three-manifolds with positive Ricci curvature,
{\sl Journal of Differential Geometry} {\bf 17} (1982), 255--306.
\sm
\ref{Ho} Chong-Wei Hong, A best constant and the Gaussian curvature, {\sl Proceedings of the 
American Mathematical Society} {\bf 97} (1986), 737--747.
\sm
\ref{J} J\"urgen Jost, Nonlinear methods in Riemannian and K\"ahlerian geometry, Birkh\"auser, 1988.
\sm
\ref{K} Erich K\"ahler, \"Uber eine bemerkenswerte Hermitesche Metrik, 
{\sl Abhandlungen aus dem Mathematischen Seminar der Universit\"at Hamburg} {\bf 9} (1933), 173--186.
\sm
\ref{Ko} S\polishl awomir Ko\polishl odziej, The complex Monge-Amp\`ere equation 
and pluripotential theory, {\sl Memoirs of the American Mathematical Society} {\bf 178} (2005), no. 840.
\sm
\ref{M1} Toshiki Mabuchi, K-energy maps integrating Futaki invariants,
{\sl T\^ohoku Mathematical Journal} {\bf 38} (1986), 575--593.
\sm
\ref{M2} \opcit, Some symplectic geometry on compact {K}\"ahler manifolds. {I},
{\sl Osaka Journal of Mathematics} {\bf 24} (1987), 227--252.
\sm
\ref{M3} \opcit, Multiplier Hermitian structures on K\"ahler manifolds, 
{\sl Nagoya Mathematical Journal} {\bf 170} (2003), 73--115.
\sm
\ref{Mo} J\"urgen Moser, A sharp form of an inequality by {N}. {T}rudinger,
{\sl Indiana University Mathematics Journal} {\bf 20} (1971), 1077--1092.
\sm
\ref{N1} \opcit, Multiplier ideal sheaves and \KE metrics of positive scalar
curvature, {\sl Annals of Mathematics} {\bf 132} (1990), 549--596.
\sm
\ref{N2}
Alan M. Nadel, On the absence of periodic points for the Ricci
curvature operator acting on the space of K\"ahler metrics, in {\it
Modern Methods in Complex Analysis: The Princeton Conference in
Honor of Gunning and Kohn} (T. Bloom et al., Eds.), Annals of
Mathematics Studies {\bf 137}, Princeton University Press, 1995, 273--282.
\sm
\ref{O} Enrico Onofri, On the positivity of the effective action in a theory of
random surfaces, {\sl Communications in Mathematical Physics} {\bf 86} (1982), 321--326.
\sm
\ref{OPS} Brad Osgood, Ralph Phillips, Peter Sarnak, Extremals of determinants of {L}aplacians,
{\sl Journal of Functional Analysis} {\bf 80} (1988), 148--211.
\sm
\ref{P} Peter Petersen, Riemannian geometry, Springer, 1998.
\sm
\ref{PSS} Duong-Hong Phong, Nata\v sa \v Se\v sum, Jacob Sturm,
Multiplier ideal sheaves and the K\"ahler-Ricci flow, preprint,
arxiv: math.DG/0611794 v2.
\sm
\ref{R1} Yanir A. Rubinstein, On iteration of the Ricci operator on the space of K\"ahler metrics, I,
manuscript, August \thhnotsosmall{14}, 2005, unpublished.
\sm
\ref{R2} \opcit, On energy functionals and the existence of \KE metrics, preprint,
arxiv: math.DG/0612440 v3.
\sm
\ref{R3} \opcit, The Ricci iteration and its applications, preprint, arxiv: 0706.2777 v3 [math.DG].
\sm
\ref{R4} \opcit, On the construction of Nadel multiplier ideal sheaves and the limiting behavior
of the Ricci flow, preprint, arxiv: 0708.1590 v2 [math.DG].
\sm
\ref{R5} \opcit, Some discretizations of geometric evolution equations and the Ricci iteration
on the space of K\"ahler metrics, II, preprint, in preparation.
\sm
\ref{Se} Stephen Semmes, Complex Monge-Amp\`ere and symplectic manifolds, 
{\sl American Journal of Mathematics} {\bf 114} (1992), 495--550.
\sm
\ref{S} Santiago R. Simanca, Heat flows for extremal \K metrics, 
{\sl Annali della Scuola Normale Superiore di Pisa} {\bf 4} (2005), 187--217.
\sm
\ref{Si} Yum-Tong Siu, Lectures on Hermitian-Einstein metrics for stable bundles
and \KE metrics,  Birkh\"auser, 1987. 
\sm
\ref{SW} Jian Song, Ben Weinkove, Energy functionals and canonical \K metrics,
{\sl Duke Mathematical Journal} {\bf 137} (2007), 159--184.
\sm
\ref{Th} Richard P. Thomas, 
Notes on GIT and symplectic reduction for bundles and varieties, 
in {\it Surveys in Differential Geometry: Essays in memory of S.-S. Chern}
(S.-T. Yau, Ed.), International Press, 2006, 221--273.
\sm
\ref{T1} Gang Tian, On \KE metrics on certain \K manifolds with $C_1(M)>0$,
{\sl Inventiones Mathematicae} {\bf 89} (1987), {225--246}.
\sm
\ref{T2} \opcit,
On Calabi's conjecture for complex surfaces with positive first Chern class,
{\sl Inventiones Mathematicae} {\bf 101} (1990), {101--172}.
\sm
\ref{T3} \opcit, On stability of the tangent bundles of Fano varieties,
{\sl International Journal of Mathematics} {\bf 3} (1992), 401--413.
\sm
\ref{T4} \opcit, The K-energy on hypersurfaces and stability, 
{\sl Communications in Analysis and Geometry} {\bf 2} (1994), 239--265.
\sm
\ref{T5} \opcit,
K\"ahler-{E}instein metrics with positive scalar curvature, 
{\sl Inventiones Mathematicae} {\bf 130} (1997), {1--37}.
\sm
\ref{T6} \opcit, Canonical Metrics in \K Geometry, Birkh\"auser, 2000.
\sm
\ref{TZ1} Gang Tian, Xiao-Hua Zhu, A nonlinear inequality of Moser-Trudinger type,
{\sl Calculus of Variations} {\bf 10} (2000), 349--354.
\sm
\ref{TZ2} \opcit, Uniqueness of K\"ahler-Ricci solitons,
{\sl Acta Mathematica} {\bf 184} (2000), 271--305.
\sm
\ref{TZ3} \opcit, A new holomorphic invariant and uniqueness of K\"ahler-Ricci solitons,
{\sl Commentarii Mathematici Helvetici} {\bf 77} (2002), 297--325.
\sm
\ref{TZ4} \opcit, Convergence of \KR flow, 
{\sl Journal of the American Mathematical Society} {\bf 20} (2007), 675--699.
\sm
\ref{TY} Gang Tian, Shing-Tung Yau, Existence of \KE metrics
on complete \K manifolds and their applications to algebraic geometry,
in {\it Mathematical Aspects of String Theory} (S.-T. Yau, Ed.),
Advanced Series in Mathematical Physics {\bf 1}, World Scientific, 1987, 574--628.
\sm
\ref{Tr} Neil S. Trudinger, On imbeddings into {O}rlicz spaces and some applications,
{\sl Journal of Mathematics and Mechanics} {\bf 17} (1967), 473--483.
\sm
\ref{WW} 
Nolan R. Wallach, Frank W. Warner, Curvature forms for 2-manifolds,
{\sl Proceedings of the American Mathematical Society} {\bf 25} (1970), 712--713.
\sm
\ref{Y} Shing-Tung Yau, On the Ricci curvature of a compact \K
manifold and the Complex \MA equation, I, {\sl Communications in Pure
and Applied Mathematics} {\bf 31} (1978), 339--411.
\sm
\ref{Z} Xiao-Hua Zhu, K\"ahler-Ricci soliton type equations on compact complex
manifolds with $C_1(M)>0$, {\sl The Journal of Geometric Analysis} {\bf 10} (2000), 759--774.

\end